\documentclass[a4paper,12pt]{amsart} 

\usepackage{amssymb}
\usepackage{euscript}
\usepackage{flafter}
\usepackage{pstricks}
\usepackage{graphicx}
\usepackage{pgfplots}

\usepackage{verbatim}

\usepackage{amsmath}

\usepackage{mathrsfs} 

\usepackage{hyperref} 

\usepackage[left=2.5cm,right=2.5cm,vmargin=2.5cm]{geometry}

\renewcommand{\theequation}{\thesection.\arabic{equation}}
\newtheorem{theorem}{Theorem}[section]
\newtheorem{lemma}[theorem]{Lemma}

\newtheorem{fact}[theorem]{Fact}

\newtheorem{claim}[theorem]{Claim}

\newcommand{\eqnsection}{
\renewcommand{\theequation}{\thesection.\arabic{equation}}
    \makeatletter
    \csname  @addtoreset\endcsname{equation}{section}
    \makeatother}
\eqnsection

\def\r{{\mathbb R}}
\def\e{{\mathbb E}}
\def\p{{\mathbb P}}
\def\P{{\bf P}}
\def\E{{\bf E}}
\def\Q{{\bf Q}}
\def\z{{\mathbb Z}}

\def\T{{\mathbb T}}
\def\G{{\mathbb G}}

\def\ee{\mathrm{e}}
\def\d{\, \mathrm{d}}

\def\mathscrHV{\mathscr{H}}
\def\mathscrHS{H^{(S)}}
\def\mathscrHx{H^{(x)}}
\def\mathscrHy{H^{(y)}}
\def\mathscrHz{H^{(z)}}

\def\HH{{\mathbb H}}

\begin{document}

\begin{abstract} The biased random walk on supercritical Galton--Watson trees is known to exhibit a multiscale phenomenon in the slow regime: the maximal displacement of the walk in the first $n$ steps is of order $(\log n)^3$, whereas the typical displacement of the walk at the $n$-th step is of order $(\log n)^2$. Our main result reveals another multiscale property of biased walks: the maximal potential energy of the biased walks is of order $(\log n)^2$ in contrast with its typical size,  which is of order $\log n$.  The proof relies on analyzing the intricate multiscale structure of the potential energy.
\end{abstract}



\title{The maximal potential energy of biased random walks on trees}


\author{Yueyun Hu}
\address{Yueyun Hu, LAGA, Universit\'e Paris XIII,  93430 Villetaneuse, France}
\email{yueyun@math.univ-paris13.fr}

\author{Zhan Shi}
\address{Zhan Shi, State Key Laboratory of Mathematical Sciences, AMSS Chinese Academy of Sciences, 100190 Beijing, China}
\email{shizhan@amss.ac.cn}

\maketitle


\section{Introduction}
   \label{s:intro}

 Let $\T$ be a supercritical Galton--Watson tree rooted at $\varnothing$. Let $\omega := (\omega(x), \, x\in \T)$ be a sequence of vectors: for each $x\in \T$, $\omega(x) := (\omega(x, \, y), \, y\in \T)$ is such that $\omega(x, \, y) \ge 0$ ($\forall y\in \T$) and that $\sum_{y\in \T} \omega(x, \, y) =1$. 

Given $\omega$, we define a random walk $(X_n, \, n\ge 0)$ on $\T$, started at $X_0 = \varnothing$, with transition probabilities given by 
$$
P_\omega \{ X_{n+1} = y \, | \, X_n =x \} = \omega(x, \, y).
$$

\noindent We assume that for each pair of vertices $x$ and $y$, $\omega(x, \, y)>0$ if and only if $y \sim x$, i.e., $y$ is either a child, or the parent, of $x$; in particular, the walk is nearest-neighbour. 

We are going to study a {\it slow regime} of the random walk $(X_n, \, n\ge 0)$. In order to observe such a slow regime, the transition probabilities $\omega(x, \, y)$ are {\it random}; i.e., given a realisation of $\omega$, we run a (conditional) Markov chain $(X_n)$. So $(X_n)$ is a randomly biased walk on the Galton--Watson tree $\T$, and can also be considered as a random walk in random environment.

 We use $\P$ to denote the law of the environment $\omega$, and $\p := \P \otimes P_\omega$ the annealed probability measure. 
 It is convenient to consider $(\omega, \T)$ as a marked tree. 
  For brevity, we say ``for almost all environment $\omega$" to mean ``for almost all $(\omega, \, \T)$".

Randomly biased walks on trees have a large literature. The model is introduced by Lyons and Pemantle~\cite{lyons-pemantle}, extending the previous model of deterministically biased walks studied in Lyons~\cite{lyons90}-\cite{lyons92}. In \cite{lyons-pemantle}, a general recurrence vs.\ transience criterion is obtained; for walks on Galton--Watson trees, the question is later also studied by Menshikov and Petritis~\cite{menshikov-petritis} and Faraud~\cite{faraud}. Ben Arous and Hammond~\cite{benarous-hammond} prove that in some sense, randomly biased walks on $\T$ are more regular than deterministically biased walks on $\T$, preventing some ``cyclic phenomena" from happening. Often motivated by results and questions in Lyons, Pemantle and Peres~\cite{lyons-pemantle-peres-ergodic} and \cite{lyons-pemantle-peres96}, the transient case has received much research attention (\cite{elie}, \cite{elie2}, \cite{elie-vitesse}, \cite{benarous-fribergh-gantert-hammond}, \cite{benarous-fribergh-sidoravicius}, \cite{benarous-hu-olla-zeitouni}), while the recurrent case has been studied in \cite{AdR17}, \cite{andreoletti-chen}, \cite{andreoletti-debs1}, \cite{AD20}, \cite{CdRH}, \cite{chen}, \cite{dR17}, \cite{faraud}, \cite{gyzbiased}, \cite{y17}, \cite{yztree},  \cite{yzslow}, \cite{K23+}, \cite{K24+} and \cite{K24++}. For a more general account of study on biased walks on trees and random walks in random environments, we refer to the books by Lyons and Peres~\cite{lyons-peres}, R\'ev\'esz \cite{revesz}, as well as the Saint-Flour lecture notes of \cite{peres} and \cite{stf}. Biased random walks in other types of random environments are also studied in the literature, such as 
\cite{AGSS23+} for biased random walks on dynamical percolation clusters, or \cite{biskup} for random walks among dynamical random conductances.

Although it is not necessary, we add a special vertex, ${\buildrel \leftarrow \over \varnothing}$, which is the parent of $\varnothing$; this simplifies our representation. The values of the transition probabilities at a finite number of vertices bringing no change to results of the paper, we can modify the value of $\omega(\varnothing, \, \bullet)$, the transition probability at $\varnothing$, in such a way that $(\omega(x, \, y), \, y\sim x)$, for $x\in \T$, are an i.i.d.\ family of random variables. 

A crucial notion in the study of the behaviour of the random walk $(X_n)$ is the {\bf potential} on $\T$, which we define by $V({\buildrel \leftarrow \over \varnothing}):=0$, $V(\varnothing):=0$ and
\begin{equation}
    V(x) 
    := 
    -
    \sum_{y\in \, ]\!] \varnothing,\, x]\!]}
    \log \,
    \frac{\omega({\buildrel \leftarrow \over y},
    \, y)}{\omega({\buildrel \leftarrow \over y}, \,
    {\buildrel \Leftarrow \over y})},
    \qquad x\in \T\backslash\{ \varnothing\},
    \label{V}
\end{equation}

\noindent where ${\buildrel \Leftarrow \over y}$ is the parent of
${\buildrel \leftarrow \over y}$, and $\, ]\!] \varnothing, \, x]\!] := [\![ \varnothing, \, x]\!] \backslash \{ \varnothing\}$, with $[\![ \varnothing, \, x]\!]$ denoting the set of vertices on the unique shortest path connecting $\varnothing$ to $x$.

Since $(\omega(x, \, y), \, y\sim x)$, for $x\in \T$, are i.i.d., the potential process $(V(x), \, x\in \T)$ is a branching random walk in the usual sense of Biggins~\cite{biggins77}, and is also well studied in the physics literature (see for example Derrida and Spohn~\cite{derrida-spohn}, represented in terms of directed polymers on trees). 

 Throughout the paper, we assume 
\begin{equation}
    \E \Big( \sum_{x: \, |x|=1} \ee^{-V(x)} \Big) =1,
    \qquad
    \E \Big( \sum_{x: \, |x|=1} V(x) \ee^{-V(x)} \Big) =0.
    \label{cond-hab}
\end{equation}

\noindent We also assume the existence of $\delta>0$ such that
\begin{equation}
    \E \Big( \sum_{x: \, |x|=1} \ee^{-(1+\delta)V(x)} \Big)
    +
    \E \Big( \sum_{x: \, |x|=1} \ee^{ \delta V(x)} \Big) 
    +
    \E \Big[ \Big( \sum_{x: \, |x|=1} 1\Big)^{1+\delta} \, \Big] <\infty .
    \label{integrability-assumption}
\end{equation}

A general result of Lyons and Pemantle~\cite{lyons-pemantle}, applied to our special setting of the Galton--Watson tree, implies that under (\ref{cond-hab}), the random walk $(X_n)$ is almost surely recurrent. This is proved in \cite{lyons-pemantle} under an additional condition on the exchangeability of 
$(V(x), \, |x|=1)$; the condition is removed in Faraud~\cite{faraud}. See also Menshikov and Petritis~\cite{menshikov-petritis} for another proof, using Mandelbrot's multiplicative cascades, modulo some additional assumptions. In the language of branching random walks, (\ref{cond-hab}) refers to the ``boundary case" in the sense of Biggins and Kyprianou~\cite{biggins-kyprianou05}. In the boundary case, the biased walk $(X_n)$ exhibits slow movement and reveals a somehow surprising multiscale behaviour: under (\ref{cond-hab}) and (\ref{integrability-assumption}) and upon the system's survival, $\max_{0\le i\le n} |X_i|$ is of order $(\log n)^3$ whereas the typical size of $|X_n|$ is of order $(\log n)^2$:  \begin{equation}
    \lim_{n\to \infty} \frac{1}{(\log n)^3} \max_{0\le i\le n} |X_i|
    =
    \frac{8}{3\pi^2 \sigma^2},
    \qquad\hbox{\rm $\p^*$-a.s.},
    \label{gyzbiased}
\end{equation}

\noindent where $\p^*:=\p(\cdot | \hbox{\rm non-extinction})$ and 
\begin{equation}
    \sigma^2
    :=
    \E \Big( \sum_{|x|=1} V(x)^2 \ee^{-V(x)} \Big) \in (0, \, \infty).
    \label{sigma}
\end{equation}

\noindent On the other hand, it is shown in \cite{yzlocaltree} that  for any $t>0$,  \begin{equation}    
 \lim_{n\to \infty} \p^* \Big( \frac{\sigma^2\, |X_n|}{(\log n)^2} \le t  \Big) 
 =
 \int_0^t \frac{1}{(2\pi s)^{1/2}} \, 
 \P \Big( {\tt m_1^{\#}} \le \frac{1}{s^{1/2}} \Big) d s\, , \label{yzlocaltree}
\end{equation} 

\noindent where ${\tt m_1^{\#}}:= \sup_{0\le r\le s \le 1} ({\tt m}_r- {\tt m}_s)$, and $({\tt m}_s)_{0\le s\le 1}$ is a Brownian meander (for the definition and basic properties of the Brownian meander, see Yen and Yor~\cite{yen-yor}). 

%
%
%
%

In dimension 1 (which corresponds heuristically to the case that every vertex has one child), a well known result of Sinai~\cite{sinai} tells that $\frac{X_n}{(\log n)^2}$ converges weakly to a non-degenerate limit; so \eqref{yzlocaltree} can be considered as a kind of companion of Sinai's theorem for the Galton--Watson tree.

Indeed, the proof of \eqref{yzlocaltree} relies on a localization of $X_n$ by some barriers constructed from the potential $V$. Using such a localization (\cite[Theorem 2.1]{yzlocaltree}), we may deduce  that assuming $(\ref{cond-hab})$ and $(\ref{integrability-assumption})$, under $\p^*$,   $\frac{V(X_n)}{\log n}$ converges in law to a  finite and  positive random variable. For any $t>0$,  under $\P ( \, \cdot \, | \, \hbox{\rm non-extinction})$,  \begin{equation} \label{cvVXn}  \lim_{n\to \infty} P_\omega\Big( \frac{V(X_n)}{\log n} \le t  \Big) 
 =\left(\frac{2}{\pi}\right)^{1/2} \, \E\left( \min\left( \frac1{\tt m_1^{\#}}, \frac{t}{\tt m_1}\right)\right) , \qquad \mbox{in probability}.   \end{equation} 

\noindent We refer to Appendix \ref{s:cvVXn} for a proof of \eqref{cvVXn}. Note that $\E ( \frac1{\tt m_1^{\#}} ) = ( \frac{\pi}{2} )^{1/2}$, see \cite{HSY-2015},   so both right-hand sides in \eqref{yzlocaltree} and \eqref{cvVXn} correspond to the partition functions of certain probability distributions on $(0, \, \infty)$. 

\medskip
In this paper, we are interested in the {\bf maximal potential energy}, 
$$
\max_{0\le k\le n} V(X_k) \, ,
$$

\noindent of the random walk $(X_i)$ in the first $n$ steps. In the literature, results on the maximal energy of random walks in random environment or related models are obtained in the one-dimensional case by Monthus and Le Doussal~\cite{monthus-ledoussal}, and for the Metropolis algorithm by Aldous~\cite{aldous}, and also by Maillard and Zeitouni~\cite{maillard-zeitouni}.

In the one-dimensional recurrent case, it is proved by Monthus and Le Doussal~\cite{monthus-ledoussal} that $\log n$ is the common order of magnitude for both $V(X_n)$ and $\max_{0\le k\le n} V(X_k)$.
 In fact, for one-dimensional random walks in random environment, it is known (Sinai~\cite{sinai}) that in $n$ steps, the maximal potential energy is bounded by $(1+o(1))\log n$ (for $n\to \infty$); more precisely, the ratio between the maximal potential energy and $\log n$ converges to a non-degenerate random variable taking values in $[0, \, 1]$. For the tree-valued random walk $(X_i)$, its restriction to each branch of $\T$ being a one-dimensional walk in random environment, the maximal potential energy along a given branch is thus bounded by $(1+o(1))\log n$, for $n\to \infty$. 

Let us present the main result of the paper.

\medskip
\begin{theorem}
\label{t:main}
 Assume $(\ref{cond-hab})$ and $(\ref{integrability-assumption})$. We have, on the set of non-extinction,
 $$
 \lim_{n\to \infty}
 \frac{1}{(\log n)^2}\, \max_{0\le k\le n} V(X_k)
 = 
 \frac12 \, ,
 \qquad\hbox{\rm $\p$-a.s.}
 $$

\end{theorem}

 \medskip
%
%
%
%
%
%



  Theorem \ref{t:main} can be viewed as another multiscale phenomenon of the randomly biased walk $(X_n)$ from the perspective of potential energy. 

The rest of the paper is as follows. Section \ref{s:preliminaries} recalls some known techniques of branching random walks which are going to be used in the proof of Theorem \ref{t:main}. The section is preceded by a brief Section \ref{s:outline_proof}, where we outline the main ideas in the proof of Theorem \ref{t:main}. It turns out that the proof relies essentially on a quenched tail estimate of excursion heights of biased walks. This tail estimate, stated in (\ref{key-estimate}), is proved in Section \ref{s:ub} by means of a second moment argument. The second moment argument being rather involving, we present it by means of two lemmas (Lemmas \ref{l:Sigma} and \ref{l:Z}), serving as the key step in the proof of the upper and lower bounds, respectively, in (\ref{key-estimate}). Lemma \ref{l:Z} is quite technical; its proof is the heart of the paper. 

Throughout the paper, we write $f(r) \sim g(r)$, $r\to \infty$, to denote $\lim_{r\to \infty} \frac{f(r)}{g(r)} =1$, and $f(r)=o(1)$, $r\to \infty$, to denote $\lim_{r\to \infty} f(r) =0$. For any pair of vertices $x$ and $y$ in $\T$, we write $x<y$ (or $y>x$) to say that $y$ is a descendant of $x$, and $x\le y$ (or $y\ge x$) to say that $y$ is either a descendant of $x$ or is $x$ itself. 

\section{Proof of Theorem \ref{t:main}: an outline}
   \label{s:outline_proof}

We assume (\ref{cond-hab}) and (\ref{integrability-assumption}), and briefly describe the proof of Theorem \ref{t:main}. Let $\varrho_0 :=0$ and let
\begin{equation}
    \varrho_n
    :=
    \inf\{ i> \varrho_{n-1}: \, X_i = {\buildrel \leftarrow \over \varnothing} \},
    \qquad n\ge 1.
    \label{rho_n}
\end{equation}

\noindent In words, $\varrho_n$ denotes the $n$-th hits to ${\buildrel \leftarrow \over \varnothing}$ by the walk $(X_i)$. It turns out that $\varrho_n = n^{1+o(1)}$ $\p$-a.s.\ for $n\to \infty$:

\begin{lemma}
\label{l:rho_n}
 Assume $(\ref{cond-hab})$ and $(\ref{integrability-assumption})$. On the set of non-extinction,
 $$
 \lim_{n\to \infty}
 \frac{\log \varrho_n}{\log n}
 = 
 1 \, ,
 \qquad\hbox{\rm $\p$-a.s.}
 $$

\end{lemma}

\medskip

The lemma is (implicitly) in \cite{yzslow} or \cite{andreoletti-debs1}. We present the proof at the end of this section, for the sake of completeness, and also to justify the passage from hitting times at $\varnothing$ to hitting times at ${\buildrel \leftarrow \over \varnothing}$. 

In view of Lemma \ref{l:rho_n}, Theorem \ref{t:main} is equivalent to the following estimate: for $\P$-almost all $\omega$ in the set of non-extinction,
\begin{equation}
    \frac{1}{(\log n)^2} \, \max_{0\le k\le \varrho_n} V(X_k)
    \to
    \frac12 \, ,
    \qquad\hbox{\rm $P_\omega$-a.s.}
    \label{proof:preparation1}
\end{equation}

At this stage, we recall an elementary result:

\begin{fact}
\label{f:valeurs-extremes}
 Let $\alpha>0$.
 Let $(\xi_n)_{n\ge 1}$ be a sequence of i.i.d.\ real-valued random variables such that $\P( \xi_1 \ge u) = \exp[ - (\alpha +o(1)) u]$, $u\to \infty$. Then 
 $$
 \lim_{n\to \infty}
 \frac{1}{\log n} \max_{1\le k\le n} \xi_k 
 = 
 \frac{1}{\alpha} \, ,
 \qquad\hbox{\rm $\P$-a.s.}
 $$

\end{fact}

\medskip

Let us go back to (\ref{proof:preparation1}). For given $\omega$, $\max_{0\le k\le \varrho_n} V(X_k)$ is the maximum of $n$ independent copies of $\max_{0\le k\le \varrho_1} V(X_k)$; so applying Fact \ref{f:valeurs-extremes} to $\xi:= [\max_{0\le k\le \varrho_1} V(X_k)]^{1/2}$ (on the set of non-extinction) and $\alpha:= 2^{1/2}$, we see that the proof of (\ref{proof:preparation1}) is reduced to verifying the following: for $\P$-almost all $\omega$ in the set of non-extinction,
\begin{equation}
    P_\omega \Big( \max_{0\le k\le \varrho_1} V(X_k) \ge r \Big)
    =
    \exp\Big( - (1+o(1)) \, (2r)^{1/2} \Big),
    \qquad r\to \infty.
    \label{proof:preparation2}
\end{equation}

For any $r>0$, let us consider the following subset of the genealogical tree:
\begin{equation}
    \mathscrHV_r
    :=
    \{ x\in \T: \, V(x) \ge r, \; \overline{V}({\buildrel \leftarrow \over x}) < r\},
    \label{M}
\end{equation}

\noindent where ${\buildrel \leftarrow \over x}$ denotes as before the parent of $x$, and for any vertex $y\in \T$, 
\begin{equation}
    \overline{V}(y)
    :=
    \max_{z\in [\![ \varnothing, \, y]\!]} V(z) ,
    \label{overline{V}}
\end{equation}

\noindent which is the maximal value of the potential $V(\cdot)$ along the path $[\![ \varnothing, \, y]\!]$. 

By definition, $\{ \max_{0\le k\le \varrho_1} V(X_k) \ge r \} = \{ T_{\mathscrHV_r} < T_{{\buildrel \leftarrow \over \varnothing}}\}$, where 
\begin{eqnarray}
    T_{\mathscrHV_r} 
 &:=& \inf\{ i\ge 0: \, X_i \in \mathscrHV_r\},
    \label{TMr}
    \\
    T_{{\buildrel \leftarrow \over \varnothing}}
 &:=& \inf\{ i\ge 0: \, X_i = {\buildrel \leftarrow \over \varnothing} \}
    =
    \varrho_1.
    \label{T_root}
\end{eqnarray}

\noindent In words, $T_{\mathscrHV_r}$ is the first hitting time of the set $\mathscrHV_r$ by the biased walk $(X_i)$. We mention that $\mathscrHV_r$ depends only on the environment, whereas $T_{\mathscrHV_r}$ involves also the behaviour of the biased walk.

We observe that (\ref{proof:preparation2}) is equivalent to saying that $\P$-almost surely on the set of non-extinction,
\begin{equation}
    P_\omega ( T_{\mathscrHV_r} < T_{{\buildrel \leftarrow \over \varnothing}})
    =
    \exp\Big( - (1+o(1)) \, (2r)^{1/2} \Big),
    \qquad r\to \infty.
    \label{key-estimate}
\end{equation}

\noindent It is (\ref{key-estimate}) we are going to prove, in Section \ref{s:ub}.

Let us close this section with the proof of Lemma \ref{l:rho_n}.

\medskip

\noindent {\it Proof of Lemma \ref{l:rho_n}.} For any $j\ge 1$, we have 
$$
P_\omega \Big\{ \max_{0\le i\le \varrho_1} |X_i| \ge j\Big\}
=
\sum_{k=1}^\infty P_\omega \Big\{ \max_{0\le i\le \varrho_1} |X_i| \ge j, \; \sum_{i=1}^{\varrho_1} {\bf 1}_{\{ X_i = \varnothing\} } = k \Big\} \, .
$$

\noindent Observe that
$$
P_\omega\Big\{ \sum_{i=1}^{\varrho_1} {\bf 1}_{\{ X_i = \varnothing\} } = k \Big\}
=
[1-\omega(\varnothing, \, {\buildrel \leftarrow \over \varnothing})]^k\,  \omega(\varnothing, \, {\buildrel \leftarrow \over \varnothing}) \, ,
$$

\noindent and that
\begin{eqnarray*}
    P_\omega \Big\{ \max_{0\le i\le \varrho_1} |X_i| \ge j \, \Big| \, \sum_{i=1}^{\varrho_1} {\bf 1}_{\{ X_i = \varnothing\} } = k\Big\} 
 &=& 1- \Big( 1- P_\omega \Big\{ \max_{0\le i\le \varrho^\varnothing} |X_i| \ge j \, \Big| \, |X_1|=1 \Big\}\, \Big)^k
    \\
 &=&1- \Big( 1- \frac{P_\omega \{ \max_{0\le i\le \varrho^\varnothing} |X_i| \ge j\} }{1-\omega(\varnothing, \, {\buildrel \leftarrow \over \varnothing})}\Big)^k \, ,
\end{eqnarray*}

\noindent where $\varrho^\varnothing := \inf\{ i\ge 1: \, X_i =\varnothing\}$. Thus
\begin{eqnarray*}
    P_\omega \Big\{ \max_{0\le i\le \varrho_1} |X_i| \ge j\Big\}
 &=& \sum_{k=1}^\infty \Big( P_\omega \{ \max_{0\le i\le \varrho^\varnothing} |X_i| \ge j\} \Big)^k \,\omega(\varnothing, \, {\buildrel \leftarrow \over \varnothing})
    \\
 &=&\frac{P_\omega \{ \max_{0\le i\le \varrho^\varnothing} |X_i| \ge j\} }{\omega(\varnothing, \, {\buildrel \leftarrow \over \varnothing}) + P_\omega \{ \max_{0\le i\le \varrho^\varnothing} |X_i| \ge j\} }
    \\
 &\le& \frac{P_\omega \{ \max_{0\le i\le \varrho^\varnothing} |X_i| \ge j\} }{\omega(\varnothing, \, {\buildrel \leftarrow \over \varnothing})} \, .
\end{eqnarray*}

\noindent So for any $n\ge 1$,
\begin{eqnarray*}
    P_\omega \Big\{ \max_{0\le i\le \varrho_n} |X_i| \ge j\Big\}
 &=&1- \Big[ 1-P_\omega \Big\{ \max_{0\le i\le \varrho_1} |X_i| \ge j\Big\} \Big]^n
    \\
 &\le& 1- \Big[ 1-\frac{P_\omega \{ \max_{0\le i\le \varrho^\varnothing} |X_i| \ge j\} }{\omega(\varnothing, \, {\buildrel \leftarrow \over \varnothing})} \Big]^n \, .
\end{eqnarray*}

\noindent By \cite{gyzbiased}, $\frac{1}{j^{1/3}} \log P_\omega \{ \max_{0\le i\le \varrho^\varnothing} |X_i| \ge j \}
\to - ( \frac{3\pi^2 \sigma^2}{8})^{1/3}$ (for $j\to \infty$) $\P$-almost surely on the set of non-extinction. Taking $j := \lceil (1+\varepsilon)^3 \frac{8}{3\pi^2 \sigma^2} (\log n)^3 \, \rceil$ with $\varepsilon>0$, we immediately see that $\P$-a.s.\ on the set of non-extinction, $\sum_\ell P_\omega \{ \max_{0\le i\le \varrho_{n_\ell} } |X_i| \ge (1+\varepsilon)^3 \frac{8}{3\pi^2 \sigma^2} (\log n_\ell)^3 \} <\infty$ if we take the subsequence $n_\ell := \lfloor \ell^{2/\varepsilon}\rfloor$, $\ell \ge 1$. By the Borel--Cantelli lemma, this yields that $\p$-almost surely, on the set of non-extinction and for all sufficiently large $\ell$,
$$
\max_{0\le i\le \varrho_{n_\ell} } |X_i|
<
(1+\varepsilon)^3 \frac{8}{3\pi^2 \sigma^2} (\log n_\ell)^3 \, ,
$$

\noindent which, in turn, implies that for $n\in [n_{\ell-1}, \, n_\ell]$,
$$
\max_{0\le i\le \varrho_n} |X_i|
<
(1+\varepsilon)^3 \frac{8}{3\pi^2 \sigma^2} (\log n_\ell)^3 
\le
(1+2\varepsilon)^3 \frac{8}{3\pi^2 \sigma^2} (\log n)^3 \, .
$$

\noindent Therefore, on the set of non-extinction,
$$
\limsup_{n\to \infty} \frac{1}{(\log n)^3} \max_{0\le i\le \varrho_n} |X_i|
\le
\frac{8}{3\pi^2 \sigma^2} \, ,
\qquad \hbox{\rm $\p$-a.s.}
$$

On the other hand, since $\varrho_n \to \infty$ $\p$-a.s., it follows from (\ref{gyzbiased}) that on the set of non-extinction,
$$
\liminf_{n\to \infty} \frac{1}{(\log \varrho_n)^3} \max_{0\le i\le \varrho_n} |X_i|
\ge
\frac{8}{3\pi^2 \sigma^2} \, ,
\qquad\hbox{\rm $\p$-a.s.}
$$

\noindent Combining the last two displayed formulas yields $\limsup_{n\to \infty} \frac{\log \varrho_n}{\log n} \le 1$ $\p$-a.s.\ on the set of non-extinction. This is the desired upper bound in Lemma \ref{l:rho_n}. The lower bound is trivial since $\varrho_n \ge 2n-1$, $\forall n\ge 1$.\hfill$\Box$

\section{Preliminaries: spinal decompositions}
   \label{s:preliminaries}
 
 We recall a useful consequence of the spinal decomposition for branching random walks. The idea of the spinal decomposition, of which we find roots in \cite{kahane-peyriere} and \cite{bingham-doney}, has been developed in the literature independently by various groups of researchers in different contexts and forms. We use here the formulation of Lyons, Pemantle and Peres~\cite{lyons-pemantle-peres} and Lyons~\cite{lyons}, based on a change-of-probabilities technique on the space of trees. We only give a brief description, referring to \cite{lyons-pemantle-peres} and \cite{lyons} for more details.

Throughout this section, we assume $\E (\sum_{|x|=1} \ee^{- V(x)} ) = 1$, which is guaranteed by (\ref{cond-hab}). Let
$$
W_n
:= 
\sum_{x: \; |x|=n} \ee^{-V(x)}, \qquad n\ge 0,
$$

\noindent which is an $(\mathscr{F}_n)$-martingale, where $\mathscr{F}_n$ denotes the $\sigma$-field generated by the branching random walk $(V(x))$ in the first $n$ generations. Kolmogorov's extension theorem ensures the existence of a probability measure $\Q$ on $\mathscr{F}_\infty$, the $\sigma$-field generated by the entire branching random walk, such that for any $n$ and any $A\in \mathscr{F}_n$,
\begin{equation}
    \Q (A)
    =
    \E( W_n \, {\bf 1}_A) \, .
    \label{Q}
\end{equation}

\noindent The distribution of $(V(x))$ under the new probability $\Q$ is called the distribution of a {\it size-biased} branching random walk. It is immediately observed that the size-biased branching random walk survives with probability one. For future use, we record here a consequence of H\"older's inequality: the assumption (\ref{integrability-assumption}) implies the existence of a constant $c_1>0$ such that
\begin{equation}
    \E_\Q \Big[ \Big( \sum_{x: \, |x|=1} \ee^{-V(x)} \Big)^{c_1} \Big]
    =
    \E \Big[ \Big( \sum_{x: \, |x|=1} \ee^{-V(x)} \Big)^{1+c_1} \Big]
    <\infty \, .
    \label{nu}
\end{equation}

We identify a branching random walk $(V(x))$ with a marked tree. On the enlarged probability space formed by marked trees with distinguished rays,\footnote{Strictly speaking, the enlarged probability space is a product space: the first coordinate concerns the branching random walk, and the second concerns the distinguished ray (= spine). In order   to keep the notation as simple as possible, we choose to work formally on the same space,  while bearing in mind that the spine $(w_n)$ is not measurable with respect to the $\sigma$-field generated by the branching random walk.} it is possible to construct a probability $\Q$ satisfying (\ref{Q}), and an infinite ray $\{ w_0= \varnothing, w_1, \dots , w_n, \dots \}$ (i.e., $w_n$ is the parent of $w_{n+1}$, and $|w_n| =n$, $\forall n\ge 0$) such that for any $n\ge 0$ and any vertex $x$ with $|x|=n$,
\begin{equation}
    \Q \{ w_n = x \, | \, \mathscr{F}_n \} 
    =
    \frac{\ee^{-V(x)}}{W_n} \, .
    \label{spine}
\end{equation}

\noindent Let us write from now on
$$
S_n
:=
V(w_n), \; n\ge 0.
$$

For any vertex $x\in \T \backslash \{ \varnothing\}$, we define
\begin{equation}
    \Delta V(x)
    :=
    V(x) - V({\buildrel \leftarrow \over x}) \, .
    \label{DeltaV}
\end{equation}

\noindent Let $f: \, \r\to [0, \, \infty)$ be a Borel function, and write 
$$
\eta_i^{(f)} 
:= 
\sum_{y: \, {\buildrel \leftarrow \over y} = w_{i-1}} f(\Delta V(y)) \, .
$$

\noindent [In particular, $\eta_1^{(f)}:= \sum_{y: \, |y|=1} f(V(y))$.] According to the spinal decomposition (see Lyons~\cite{lyons}), $(S_i-S_{i-1}, \, \eta_i^{(f)})$, $i\ge 1$, are i.i.d.\ under $\Q$. 

For any vertex $x\in \T$, let $x_i$ be the ancestor of $x$ in the $i$-th generation for $0\le i\le |x|$ (so $x_0=\varnothing$, and $x_{|x|}=x$). Let $n\ge 1$, and let $g: \, \r^{2n} \to [0, \, \infty)$ be a Borel function. By definition of $\Q$, we have
\begin{eqnarray*}
 &&\E \Big[ \sum_{x: \, |x|=n} g\Big( V(x_i), \; \sum_{y: \, {\buildrel \leftarrow \over y} = x_{i-1}} f(\Delta V(y)), \; 1\le i\le n\Big) \Big]
    \\
 &=&\E_\Q \Big[ \frac{1}{W_n} \sum_{x: \, |x|=n} g\Big( V(x_i), \; \sum_{y: \, {\buildrel \leftarrow \over y} = x_{i-1}} f(\Delta V(y)), \; 1\le i\le n\Big) \Big] ,
\end{eqnarray*}

\noindent which, according to (\ref{spine}), is 
\begin{eqnarray*}
 &=& \E_\Q \Big[ \sum_{x: \, |x|=n} \ee^{V(x)} \, {\bf 1}_{\{ w_n=x\} } \,
    g\Big( V(x_i), \; \sum_{y: \, {\buildrel \leftarrow \over y} = x_{i-1}} f(\Delta V(y)), \; 1\le i\le n\Big) 
    \Big] 
    \\
 &=&\E_\Q \Big[ \ee^{V(w_n)} \, g\Big( V(w_i), \; \sum_{y: \, {\buildrel \leftarrow \over y} = w_{i-1}} f(\Delta V(y)), \; 1\le i\le n\Big) \Big] .
\end{eqnarray*}

\noindent In our notation, this means that
\begin{eqnarray}
 &&\E \Big[ \sum_{x: \, |x|=n} g\Big( V(x_i), \; \sum_{y: \, {\buildrel \leftarrow \over y} = x_{i-1}} f(\Delta V(y)), \; 1\le i\le n\Big) \Big]
    \nonumber
    \\
 &=& \E_\Q \Big[ \ee^{S_n} \, g\Big( S_i, \; \eta_i^{(f)}, \; 1\le i\le n\Big) \Big] .
    \label{spinal-decomposition}
\end{eqnarray}

\noindent A special case of (\ref{spinal-decomposition}) is of particular interest: for any $n\ge 1$ and any Borel function $g: \r^n \to \r_+$,
\begin{equation}
    \E\Big[ \sum_{x: \, |x|=n} g(V(x_1), \, \cdots, \, V(x_n)) \Big]
    =
    \E_{\Q}\Big[ \ee^{S_n} \, g(S_1, \, \cdots, \, S_n) \Big] .
    \label{many-to-one}
\end{equation}

\noindent This is the so-called many-to-one formula, and can also be directly checked by induction on $n$ without using (\ref{spine}). An immediate consequence of (\ref{many-to-one}) is that the assumption (\ref{cond-hab}) yields $\E_{\Q}(S_1)=0$, whereas the assumption (\ref{integrability-assumption}) implies that
$$
    \E_{\Q}(\ee^{a\, S_1})
    <
    \infty,
    \qquad \forall 0\le a < \delta \, .
$$

\noindent The existence of some finite exponential moments allows us to use the last displayed formula\footnote{More precisely, we apply the formula of Chang~\cite{chang} to the ladder height of our mean-zero random walk via the Theorem on page 250 of Doney~\cite{doney}.} on page 1229 of Chang~\cite{chang} to see that there exists a constant $c_2>0$ satisfying
\begin{equation}
    \sup_{b>0} \E_{\Q} \Big[ \exp( c_2\, \Delta S_{\mathscrHS_b}) \Big]
    <
    \infty \, ,
    \label{moment-expo-overshoot}
\end{equation}

\noindent where
\begin{eqnarray}
    \Delta S_i 
 &:=& S_i - S_{i-1},
    \qquad i\ge 1 ,
    \label{DeltaS}
    \\
    \mathscrHS_r
 &:=& \inf\{ i\ge 0: \, S_i \ge r\} ,
    \qquad r\ge 0 \, .
    \label{mathscrHS_r}
\end{eqnarray}

The formula (\ref{spinal-decomposition}), stated for any given generation $n$, remains valid if $n$ is replaced by $\mathscrHV_r$, with $\mathscrHV_r := \{ x\in \T: \, V(x) \ge r, \; \overline{V}({\buildrel \leftarrow \over x}) < r\}$ as in (\ref{M}). Indeed, according to Proposition 3 of \cite{AHZ}, for any $r>0$ and any measurable functions $f$ and $g$,
\begin{eqnarray}
 &&\E \Big[ \sum_{x\in \mathscrHV_r} g\Big( V(x_i), \; \sum_{y: \, {\buildrel \leftarrow \over y} = x_{i-1}} f(\Delta V(y)), \; 1\le i\le |x| \Big) \Big]
    \nonumber
    \\
 &=& \E_\Q \Big[ \exp(S_{\mathscrHS_r})\, g\Big( S_i, \; \eta_i^{(f)} , \; 1\le i\le \mathscrHS_r\Big) \Big] ,
    \label{decomposition-epinale:stopping-lines}
\end{eqnarray}

\noindent where $\eta_i^{(f)} := \sum_{y: \, {\buildrel \leftarrow \over y} = w_{i-1}} f(\Delta V(y))$ as before. We recall that $(S_i-S_{i-1}, \, \eta_i^{(f)})$, $i\ge 1$, are i.i.d.\ under $\Q$.

In particular, we have the following analogue of the many-to-one formula for $\mathscrHV_r$:
\begin{equation}
    \E\Big[ \sum_{x\in \mathscrHV_r} g(V(x_1), \, \cdots, \, V(x_{|x|})) \Big]
    =
    \E_{\Q} \Big[ \exp(S_{\mathscrHS_r}) \, g(S_1, \, \cdots, \, S_{\mathscrHS_r}) \Big] \, .
    \label{many-to-one:stopping-lines}
\end{equation}

\section{Multiscale potential energy optimization and proof of Theorem \ref{t:main}}
   \label{s:ub}

The proof of Theorem \ref{t:main} relies on a multiscale optimization analysis of the potential theory. The main ingredients are summarized in a couple of lemmas, stated as Lemmas \ref{l:Sigma} and \ref{l:Z} below. Lemma \ref{l:Z}, rather technical, consists of three estimates, namely, (\ref{E(Z)}), (\ref{E(Z2)}) and (\ref{(EZ)2}). Here is how the proofs are organized:

\medskip

\qquad$\bullet$ Subsection \ref{subs:proof-thm}: proof of Theorem \ref{t:main}, by admitting Lemmas \ref{l:Sigma} and \ref{l:Z}.

\qquad$\bullet$ Subsection \ref{subs:proof-first-lemma}: proof of Lemma \ref{l:Sigma}.

\qquad$\bullet$ Subsection \ref{subs:proof-second-lemma-(i)}: proof of Lemma \ref{l:Z}, the inequality (\ref{E(Z)}).

\qquad$\bullet$ Subsection \ref{subs:proof-second-lemma-(ii)}: proof of Lemma \ref{l:Z}, the inequality (\ref{E(Z2)}).

\qquad$\bullet$ Subsection \ref{subs:proof-second-lemma-(iii)}: proof of Lemma \ref{l:Z}, the inequality (\ref{(EZ)2}).

\medskip

Throughout the section, we assume $(\ref{cond-hab})$ and $(\ref{integrability-assumption})$. 

For any $x\in \T \cup \{ {\buildrel \leftarrow \over \varnothing} \}$, let
\begin{equation}
    T_x 
    := 
    \inf\{ n\ge 0: \, X_n =x\} ,
    \qquad (\inf\varnothing := \infty)
    \label{Tx}
\end{equation}

\noindent which stands for the first hitting time of the vertex $x$ by the biased walk. [In the special case $x:= {\buildrel \leftarrow \over \varnothing}$, (\ref{Tx}) is in agreement with (\ref{T_root}).] For $r>0$, we recall from (\ref{TMr}) that
$$
T_{\mathscrHV_r} 
:=
\inf\{ i\ge 0: \, X_i \in \mathscrHV_r\} ,
$$

\noindent where $\mathscrHV_r := \{ x\in \T: \, V(x) \ge r, \; \overline{V}({\buildrel \leftarrow \over x}) < r\}$ as in (\ref{M}). 

Our first preliminary result is as follows. Its innocent-looking statement masks the multiscale nature.  


\begin{lemma}
\label{l:Sigma}

 Assume $(\ref{cond-hab})$ and $(\ref{integrability-assumption})$. We have\footnote{Of course, $\E [P_\omega(\cdots)]$ is nothing else but $\e(\cdots)$.}
 $$
 \limsup_{r\to \infty} \frac{1}{(2r)^{1/2}} \, \log \E [P_\omega (T_{\mathscrHV_r} < T_{{\buildrel \leftarrow \over \varnothing}} )]
 \le
 -1 \, .
 $$

\end{lemma}

\medskip

The statement of the second lemma requires the introduction of the basic multiscale setting in our analysis of the potential energy. Let $r>0$. Let $\chi \in (\frac12, \, 1)$. Let
\begin{eqnarray}
    k 
 &:=& \lfloor r^{1-\chi} \rfloor \, ,
    \label{k}
    \\
    h_m
 &:=& \frac{r}{k} \, m \, , 
    \qquad 0\le m\le k \, ,
    \label{hm}
    \\
    \lambda_m
 &:=& (2r)^{1/2} \, (\frac{k-m+1}{k})^{1/2}\, ,
    \quad
    1\le m\le k \, .
    \label{lambda_m}
\end{eqnarray}

\noindent [The value of $\lambda_m$ defined in \eqref{lambda_m} represents, in fact, the solution to an optimization problem.] For any $x\in \T$ and any $0\le s \le \overline{V}(x)$ (for definition of $\overline{V}(x)$, see (\ref{overline{V}})), let
\begin{equation}
    \mathscrHx_s
    =
    \inf\Big\{ i\ge 0: \, V(x_i) \ge s, \; V(x_j) < s, \; \forall j\in [0, \, i)\Big\} \, .
    \label{H(x)_r} 
\end{equation}

\noindent In words, $\mathscrHx_s$ is the generation of the oldest vertex in the path $[\![ \varnothing, \, x]\!]$ such that the value of the branching random walk $V(\cdot)$ is at least $s$.

For $x\in \mathscrHV_r:= \{ x\in \T: \, V(x) \ge r, \; \overline{V}({\buildrel \leftarrow \over x}) < r\}$, we set\footnote{As such, $a_i^{(x)}$ is well defined for all $0\le i < \mathscrHx_r = |x|$ (for $x\in \mathscrHV_r$). The value of $a_i^{(x)}$ for $i=\mathscrHx_r$ plays no role. [One can, for example, set $a_i^{(x)} := a_{i-1}^{(x)}$ for $i=\mathscrHS_r$.]}
\begin{equation}
    a_i^{(x)}
    :=
    \lambda_m,
    \qquad \hbox{ \rm if } \mathscrHx_{h_{m-1}} \le i< \mathscrHx_{h_m} \hbox{ \rm for } m\in [1, \, k] \, .
    \label{aV}
\end{equation}

In order to control the increments of the potential along the children of vertices in the spine, we introduce, for any vertex $x\in \T$, the following quantity
\begin{equation}
    \Lambda(x)
    :=
    \sum_{y: \, {\buildrel \leftarrow \over y} = x} \ee^{-\Delta V(y)}
    =
    \sum_{y: \, {\buildrel \leftarrow \over y} = x} \ee^{-[V(y)-V(x)]} \, .
    \label{Lambda} 
\end{equation}

Let $c_1>0$ be the constant in (\ref{nu}). Fix $\varepsilon> 0$, $\beta \ge 0$, $0<\varepsilon_1< c_1 \, \varepsilon$ and $\theta \in (\frac12 , \, \chi)$.\footnote{For Lemma \ref{l:Z}, we can take any $\theta \in (\frac{\chi}{2} , \, \chi)$, but condition $\max_{1\le m\le k} \Delta V(x_{\mathscrHx_{h_{m-1}}}) \le r^\theta$ is also exploited in Section~\ref{subs:proof-first-lemma} in the proof of Lemma \ref{l:Sigma}, where $\theta$ needs to be greater than $\frac12$. In order to avoid any possibility of confusion, we take $\theta \in (\frac12 , \, \chi)$ once for all.} We consider the following subset of $\mathscrHV_r$:
\begin{eqnarray}
    \mathscrHV_r^*
 &:=& \Big\{ x\in \mathscrHV_r: \, \max_{1\le m<k} \Delta V(x_{\mathscrHx_{h_m}}) \le r^\theta , \; \underline{V}(x) \ge - \beta , \; |x| < \lfloor \ee^{\varepsilon_1 \, r^{1/2}}\rfloor \, , \;
    \nonumber
    \\
 && \overline{V}(x_j) - V(x_j) \le a_j^{(x)} , \; \forall 0\le j<|x| , \;  \max_{0\le j<|x|}\Lambda(x_j) \le \ee^{\varepsilon r^{1/2}} \Big\} ,
    \label{small-overshoots}
\end{eqnarray}

\noindent where $\Delta V(y) := V(y) - V({\buildrel \leftarrow \over y})$ as in (\ref{DeltaV}), $\Lambda(x) := \sum_{y: \, {\buildrel \leftarrow \over y} = x} \ee^{-\Delta V(y)}$ as in (\ref{Lambda}), and 
$$
\underline{V}(y)
:=
\min_{z\in [\![ \varnothing, \, y]\!]} V(z) \, ,
$$

\noindent for all $y\in \T$. See Figure 1. 
\begin{figure}[h]\label{figure1}
\centering
\begin{tikzpicture}
    \begin{axis}[
       	axis y line=left,
    	axis x line=middle,
        ymin=-1, ymax=8,
        xmin=0, xmax=13,
        xtick=\empty, ytick=\empty,
        height=8cm,
        width=14cm,
        clip=false,
        ]

        \addplot[thick, blue, no markers] coordinates {
            (0,0) (0.5,-0.5) (1,0.7) (1.5, -0.3) (2,1) (2.5,1.3) 
            (3,1) (3.5,1.3) (4,2.5) (4.5,1.6) (5,3.0) (5.5,2.2) 
            (6,3.8) (6.5,3.5) (7,4.5) (7.5,4.2)
            (8,3) (8.5, 4.5)  (9,5.5) (9.5, 5.2) (10,6) (10.5, 5.6) (11,6.5) (12,7.3)
        };

        \addplot[dotted] coordinates {(0,2) (12,2)};
        \addplot[dotted] coordinates {(0,5) (12,5)};
        \addplot[dotted] coordinates {(0,7) (12,7)};
        \addplot[dotted] coordinates {(8,3) (12.3,3)};
        \addplot[dotted] coordinates {(7,4.5) (12.3,4.5)};

        \addplot[dotted] coordinates {(4,0) (4,2.3)};
        \addplot[dotted] coordinates {(9,0) (9,5.5)};
        \addplot[dotted] coordinates {(12,0) (12,7.3)};
        \addplot[<->] coordinates {(12.3,4.5) (12.3,3)};
        \node at (axis cs:-0.5,2) {\scriptsize $h_{m-1}$};
        \node at (axis cs:-0.5,5) {\scriptsize $h_{m}$};
        \node at (axis cs:-0.5,7) {\scriptsize $r$};
        \node at (axis cs:4.2,-0.4) {\scriptsize $H^{(x)}_{h_{m-1}}$};
        \node at (axis cs:9.1,-0.4) {\scriptsize $H^{(x)}_{h_{m}}$};
        \node at (axis cs:12.1,-0.4) {\scriptsize $H^{(x)}_{r}$};
         \node at (axis cs:13,-0.4) {\scriptsize $i$};
          \node at (axis cs:-0.5, 8) {\scriptsize $V(x_i)$};
        \node at (axis cs:12.8,3.8) {\scriptsize $\leq \lambda_m$};
        
        \node[anchor=east] at (axis cs:0,0) {\scriptsize 0};
    \end{axis}
\end{tikzpicture}
 \caption{\scriptsize An illustration of $(V(x_i), 0\le i \le H_r^{(x)})$ when $x\in  \mathscrHV_r^*$. We have $h_k=r$ and $1\le m\le k$.  For $H^{(x)}_{h_{m-1}} \leq i < H^{(x)}_{h_m}$, $a^{(x)}_i = \lambda_m$, and  $\overline{V}(x_i) - V(x_i) \leq \lambda_m$.}
\end{figure}
 
%
%
%
Define $Z_r=Z_r(\varepsilon, \, \varepsilon_1 , \, \beta, \, \theta, \, \chi)$ by 
\begin{equation}
    Z_r
    :=
    \sum_{x\in \mathscrHV_r^*} 
    {\bf 1}_{\{ T_x < T_{{\buildrel \leftarrow \over \varnothing}} \} } \, .
    \label{Z}
\end{equation}

\noindent The reason for which we are interested in $Z_r$ is the obvious relation $\{ T_{\mathscrHV_r} < T_{{\buildrel \leftarrow \over \varnothing}} \} \supset \{ Z_r >0\}$. 

In the definition of $Z_r$, everything depends only on the random potential $V(\cdot)$, except for $T_x$ and $T_{{\buildrel \leftarrow \over \varnothing}}$, both of which depend also on the movement of the biased random walk $(X_i)$. 

The following lemma, which is the main technical result of the paper, summarizes the moment properties of $Z_r$ we need in the multiscale analysis of the potential energy.

\begin{lemma}
\label{l:Z}

 Assume $(\ref{cond-hab})$ and $(\ref{integrability-assumption})$. 
 For any $0<\varepsilon_1<c_1 \, \varepsilon$, $\beta \ge 0$ and $\frac12 < \theta < \chi < 1$, we have
 \begin{eqnarray}
     \liminf_{r\to \infty} \frac{1}{(2r)^{1/2}} \, \log \E [ E_\omega (Z_r)]
  &\ge& 
     -1 - \frac{\varepsilon_1}{2^{1/2}}\, ,
     \label{E(Z)}
     \\
     \limsup_{r\to \infty} \frac{1}{(2r)^{1/2}} \, \log \E[ E_\omega  (Z_r^2)]
  &\le& -1 + 2^{1/2} \, (\varepsilon+\varepsilon_1) \, ,
     \label{E(Z2)}
     \\
     \limsup_{r\to \infty} \frac{1}{(2r)^{1/2}} \, \log \E[ ( E_\omega  Z_r)^2]
  &\le& -2 + 2^{1/2} \, \varepsilon \, .
     \label{(EZ)2}
 \end{eqnarray}
\end{lemma}

\medskip

By admitting Lemmas \ref{l:Sigma} and \ref{l:Z} for the time being, we are ready to prove Theorem \ref{t:main}.

\subsection{Proof of Theorem \ref{t:main}}
\label{subs:proof-thm}

We have seen in Section \ref{s:outline_proof} that the proof of Theorem \ref{t:main} consists of verifying (\ref{key-estimate}), of which we recall the statement: under the assumptions $(\ref{cond-hab})$ and $(\ref{integrability-assumption})$, $\P$-almost surely on the set of non-extinction,
$$
    \lim_{r\to \infty} \frac{1}{(2r)^{1/2}} \, \log P_\omega ( T_{\mathscrHV_r} < T_{{\buildrel \leftarrow \over \varnothing}})
    =
    -1 \, .
    \leqno(\ref{key-estimate})
$$

\noindent Lemma \ref{l:Sigma} is useful in the proof of the upper bound in (\ref{key-estimate}), and Lemma \ref{l:Z} the lower bound.

We start with the proof of the upper bound, by means of Lemma \ref{l:Sigma}. 
Let 
$$
\P^* (\, \cdot \, ) 
:= 
\P ( \, \cdot \, | \, \hbox{\rm non-extinction}) \, .
$$

\noindent By Lemma \ref{l:Sigma} and the Markov inequality, 
$$
\P^* \{ P_\omega ( T_{\mathscrHV_r} < T_{{\buildrel \leftarrow \over \varnothing}} ) 
> \ee^{-(1-\varepsilon) (2r)^{1/2}} \} 
\le 
\ee^{-c_3 \, (2r)^{1/2}} \, ,
$$

\noindent for some $c_3 = c_3(\varepsilon)>0$ and all sufficiently large $r$. An application of the Borel--Cantelli lemma 
yields that with $\P^*$-probability 
1, for all sufficiently large {\it integer} numbers $r>0$, $P_\omega ( T_{\mathscrHV_r} < T_{{\buildrel \leftarrow \over \varnothing}} ) \le \ee^{-(1-\varepsilon) (2r)^{1/2}}$. Since $r\to T_{\mathscrHV_r}$ is non-decreasing, we can remove the condition that $r$ be integer. As a consequence, 
$$
\limsup_{r\to \infty} \frac{1}{(2r)^{1/2}} \, \log P_\omega ( T_{\mathscrHV_r} < T_{{\buildrel \leftarrow \over \varnothing}})
\le
-1 \, ,
\qquad\hbox{\rm $\P^*$-a.s.,}
$$

\noindent which is the desired upper bound in (\ref{key-estimate}).

We now turn to the proof of the lower bound. 
Since 
$\E[P_\omega \{ Z_r> 0 \}] = (\P \otimes P_\omega)\{ Z_r> 0 \}$, it follows from the Cauchy--Schwarz inequality that
$$
\E[P_\omega \{ Z_r> 0 \}
]
\ge 
\frac{\{ \E [ E_\omega (Z_r)
]\}^2}{\E[ E_\omega  (Z_r^2)
]} \, .
$$

\noindent Applying (\ref{E(Z)}) and (\ref{E(Z2)}) of Lemma \ref{l:Z} yields that
\begin{equation}
    \liminf_{r\to \infty} \frac{1}{(2r)^{1/2}} \, \log \E[P_\omega \{ Z_r> 0 \}
    ]
    \ge
    -1 - 2^{1/2}\, (\varepsilon + \varepsilon_1) - 2^{1/2} \, \varepsilon_1 \, .
    \label{Cauchy-Schwarz-Z}
\end{equation}

\noindent On the other hand, by the Markov inequality, $P_\omega \{ Z_r> 0 \} \le E_\omega(Z_r)$, 
so it follows from (\ref{(EZ)2}) of Lemma \ref{l:Z} that
\begin{equation}
    \limsup_{r\to \infty} \frac{1}{(2r)^{1/2}} \, \log \E[(P_\omega \{ Z_r> 0 \}
    )^2 ]
    \le
    -2 + 2^{1/2} \, \varepsilon \, .
    \label{Markov-Z}
\end{equation}

Recall (a special case of) the Paley--Zygmund inequality: for any non-negative random variable $\xi$, we have $\P\{ \xi > \frac12 \E(\xi)\} \ge \frac14 \frac{[\E(\xi)]^2}{\E(\xi^2)}$. We apply it to $\xi := P_\omega \{ Z_r> 0 \}
$. 
In view of (\ref{Cauchy-Schwarz-Z}) and (\ref{Markov-Z}), we obtain: for any $\varepsilon_2  > 6\varepsilon + 8\varepsilon_1$ and all sufficiently large $r$,
$$
    \P \{ P_\omega \{ Z_r> 0 \}
    > \ee^{-(1+\varepsilon_2)(2r)^{1/2}} \}
    \ge
    \ee^{-\varepsilon_2\, r^{1/2}}.
$$

Let
\begin{equation}
    \gamma_r
    :=
    P_\omega ( T_{\mathscrHV_r} < T_{{\buildrel \leftarrow \over \varnothing}}).
    \label{gamma}
\end{equation}

\noindent Since $\{ T_{\mathscrHV_r} < T_{{\buildrel \leftarrow \over \varnothing}} \} \supset \{ Z_r >0\}$, we have $\gamma_r \ge P_\omega \{ Z_r> 0 \} 
$. Consequently, for all sufficiently large $r>0$,
\begin{equation}
    \P \{ \gamma_r > \ee^{-(1+\varepsilon_2)(2r)^{1/2}} \}
    \ge
    \ee^{-\varepsilon_2 \, r^{1/2}}.
    \label{Paley-Zygmund-csq}
\end{equation}

As this stage, it is convenient to have the following preliminary estimate. Recall from (\ref{overline{V}}) that $\overline{V}(x) := \max_{z\in [\![ \varnothing, \, x]\!]} V(z)$.

\begin{claim}
\label{claim:GW_mean}

 Let $c_4>0$ be a constant satisfying $(\ref{C})$ below. Let $0<\alpha<\frac12$. Let
 $$
 \mu_L
 :=
 \E \Big( \sum_{x: \, |x| = L} {\bf 1}_{\{ V(x) \ge L^\alpha\} } \, {\bf 1}_{\{ \overline{V} (x) < 2L^\alpha\} } \, {\bf 1}_{\{ \prod_{j=0}^{L-1} [1+\Lambda(x_j)] \le \ee^{c_4 L}\} } \Big) ,
 $$
 
 \noindent where $\Lambda(x) := \sum_{y: \, {\buildrel \leftarrow \over y} = x} \ee^{-\Delta V(y)}$ as in $(\ref{Lambda})$. Then $\lim_{L\to \infty} \mu_L =\infty$.
 
\end{claim}

\medskip

\noindent {\it Proof of Claim \ref{claim:GW_mean}.} By (\ref{spinal-decomposition}), we have
$$
\mu_L
=
\E_{\Q} \Big( \ee^{S_L} \, {\bf 1}_{\{ S_L \ge L^\alpha\} } \, {\bf 1}_{\{ \overline{S}_L < 2L^\alpha\} } \, {\bf 1}_{\{ \prod_{j=1}^L (1+ \eta_j) \le \ee^{c_4L}\} } \Big) ,
$$

\noindent where $(S_j-S_{j-1}, \, \eta_j)$, $j\ge 1$, are i.i.d.\ random vectors under $\Q$, with $\eta_1 := \sum_{y: \, |y|=1} \ee^{- V(y)}$, and
\begin{equation}
    \overline{S}_j
    :=
    \max_{0\le i\le j} S_i ,
    \qquad j\ge 0 \, .
    \label{S_bar}
\end{equation}

\noindent Hence
\begin{eqnarray}
    \mu_L
 &\ge&\ee^{L^\alpha} \Q \Big\{ S_L \ge L^\alpha , \, \overline{S}_L < 2L^\alpha, \, \prod_{j=1}^L (1+ \eta_j) \le \ee^{c_4L} \Big\} 
    \nonumber
    \\
 &\ge&\ee^{L^\alpha} \Big[ \Q \{ S_L \ge L^\alpha , \, \overline{S}_L < 2L^\alpha \} - \Q \Big\{ \prod_{j=1}^L (1+ \eta_j) > \ee^{c_4L} \Big\} \Big] .
    \label{mu_L>}
\end{eqnarray}

\noindent We claim that for some constants $c_5>0$ and $c_6>0$, 
\begin{eqnarray}
    \liminf_{L\to \infty} \, L^{\frac32-2\alpha}\, \Q \{ S_L \ge L^\alpha , \, \overline{S}_L < 2L^\alpha \} 
 &\ge& c_5 \, ,
    \label{mu_L>(1)} 
    \\
    \limsup_{L\to \infty} \, \frac{1}{L} \, \log \Q \Big\{ \prod_{j=1}^L (1+ \eta_j) > \ee^{c_4L} \Big\}
 &\le& - c_6 \, . 
    \label{mu_L>(2)} 
\end{eqnarray}

\noindent It is clear that Claim \ref{claim:GW_mean} will follow from (\ref{mu_L>(1)}) and (\ref{mu_L>(2)}).

To check (\ref{mu_L>(1)}), we use $\Q \{ S_L \ge L^\alpha , \, \overline{S}_L < 2L^\alpha \} \ge \Q \{ L^\alpha \le S_L < 2L^\alpha , \, \overline{S}_{L-1} \le S_L \}$. Since $(S_L-S_{L-i}, \, 0\le i\le L)$ is distributed as $(S_i, \, 0\le i\le L)$, the latter probability equals $\Q \{ L^\alpha \le S_L < 2L^\alpha , \, S_i \ge 0, \, \forall 1\le i\le L \}$, which can be written as $\Q\{ S_i \ge 0, \, \forall 1\le i\le L\} \times \Q \{ L^\alpha \le S_L < 2L^\alpha \, | \, S_i \ge 0, \, \forall 1\le i\le L \}$. It is known (Kozlov~\cite{kozlov}) that $L^{1/2} \, \Q\{ S_i \ge 0, \, \forall 1\le i\le L\}$ converges (when $L\to \infty$) to a positive limit, whereas according to Caravenna~\cite{caravenna}, $\liminf_{L\to \infty} L^{1-2\alpha} \, \Q \{ L^\alpha \le S_L < 2L^\alpha \, | \, S_i \ge 0, \, \forall 1\le i\le L \} >0$. This yields (\ref{mu_L>(1)}).

The proof of (\ref{mu_L>(2)}) is also elementary. Let $\delta_1 \in (0, \, 1]$. By the Markov inequality,
$$
\Q \Big\{ \prod_{j=1}^L (1+ \eta_j) > \ee^{c_4L} \Big\}
\le
\Big\{ \ee^{-\delta_1 c_4}\, \E_{\Q} [(1+\eta_1)^{\delta_1}] \Big\}^L
\le
\Big\{ \ee^{-\delta_1 c_4}\, [1+\E_{\Q} (\eta_1^{\delta_1})] \Big\}^L .
$$

\noindent Note that $\E_{\Q} (\eta_1^{\delta_1}) = \E_\Q [(\sum_{|y|=1} \ee^{- V(y)})^{\delta_1}] <\infty$ if we choose $\delta_1 := \min\{ c_1, \, 1\}$ (see (\ref{nu})). So, as long as 
\begin{equation}
    c_4
    > 
    \frac{\log [1+\E_{\Q} (\eta_1^{\delta_1})]}{\delta_1} \, ,
    \label{C}
\end{equation}

\noindent we have $\ee^{-\delta_1 c_4} [1+\E_{\Q} (\eta_1^{\delta_1})] <1$, which yields (\ref{mu_L>(2)}). Claim \ref{claim:GW_mean} is proved.\hfill$\Box$

\medskip

We continue with our proof of Theorem \ref{t:main}, or more precisely, of the lower bound in (\ref{key-estimate}). By Claim \ref{claim:GW_mean}, we are entitled to choose and {\bf fix} an integer $L$ such that $\mu_L>1$.

Let us construct a super-critical Galton--Watson $\G^{(L)}$ which is a sub-tree of $\T$. The vertices in $\G_1^{(L)}$, the first generation of $\G^{(L)}$, are those $x\in \T$ with $|x| = L$ such that 
$$
V(x) \ge L^\alpha \, ,
\qquad
\overline{V}(x) < 2L^\alpha \, ,
\qquad
\prod_{j=0}^{L-1} [1+ \Lambda(x_j)] \le \ee^{c_4L} \, ,
$$

\noindent where $\Lambda(x) := \sum_{y: \, {\buildrel \leftarrow \over y} = x} \ee^{-\Delta V(y)}$ as in (\ref{Lambda}). More generally, for any $n\ge 2$, the vertices in $\G_n^{(L)}$, the $n$-th generation of $\G^{(L)}$, are those $x\in \T$ with $|x| = nL$ such that $V(x) - V(x^*) \ge L^\alpha$, that $\max_{(n-1)L \le i \le nL} [V(x_i) - V(x^*)] < 2L^\alpha$ and that $\prod_{j=(n-1)L}^{nL-1} [1+\Lambda(x_j)]\le \ee^{c_4L}$, where $x^*$ is the parent in $\G_{n-1}^{(L)}$ of $x$ (so $x^* = x_{(n-1)L}$ as a matter of fact). 


Let $c_4>0$ be a constant satisfying (\ref{C}). Let $\mathscrHV_s := \{ x\in \T: \, V(x) \ge s, \; \overline{V}({\buildrel \leftarrow \over x}) < s\}$ as defined in (\ref{M}). Let
$$
\mathscr{K}_s
:=
\Big\{ x\in \mathscrHV_s: \, \prod_{j=0}^{|x|-1} [1+\Lambda(x_j)] \le \ee^{2c_4 L^{1-\alpha} s}, \; |x| \le 2L^{1-\alpha}s, \; V(x) \le 4s \Big\} \, .
$$

\noindent We need an elementary result.

\begin{claim}
\label{claim:GW_embedded}

 For $n\ge 1$ and $s\in [2nL^\alpha, \, 2(n+1)L^\alpha]$,
 \begin{equation}
     \# \mathscr{K}_s
     \ge
     \sum_{y\in \G^{(L)}_n} {\bf 1}_{\{ \exists z\in \G^{(L)}_{2n+2}: \; y<z\} }  \, .
     \label{GW-inclusion}
 \end{equation}

\end{claim}

\noindent {\it Proof of Claim \ref{claim:GW_embedded}.} Let $y\in \G^{(L)}_n$ be such that there exists $z\in \G^{(L)}_{2n+2}$ with $y<z$. By definition of $\G^{(L)}$, we have $V(y) < 2nL^\alpha \le s$ and $V(z) \ge (2n+2)L^\alpha \ge s$. So there exists $x\in [\![ y, \, z]\!]$ such that $x\in \mathscrHV_s$. Since $x$ is a descendant of $y$, all we need is to check that $x\in \mathscr{K}_s$.

Since $z\in \G^{(L)}_{2n+2}$, we have, by definition of $\G^{(L)}$, $\prod_{j=0}^{|z|-1}[1+\Lambda(z_j)] \le \ee^{c_4(2n+2)L}$, and a fortiori (using $x\le z$), $\prod_{j=0}^{|x|-1} [1+\Lambda(x_j)] \le \ee^{c_4(2n+2)L} \le \ee^{4c_4nL} \le \ee^{2c_4L^{1-\alpha}s}$.

On the other hand, $|x| \le |z| = (2n+2)L \le 4nL \le 2L^{1-\alpha}s$.

Finally, $V(x) \le (2n+2)2L^\alpha \le 8nL^\alpha \le 4s$. As a conclusion, $x\in \mathscr{K}_s$.\hfill$\Box$

\bigskip

We come back to the proof of the lower bound in (\ref{key-estimate}). We use the trivial inequality
$$
\sum_{y\in \G^{(L)}_n} {\bf 1}_{\{ \exists z\in \G^{(L)}_{2n+2}: \; y<z\} }
\ge
\sum_{y\in \G^{(L)}_n} {\bf 1}_{\{ \textrm{the sub-tree in $\G^{(L)}$ rooted at $y$ survives} \} } \, .
$$

\noindent Since $\G^{(L)}$ is supercritical, there exist constants $c_7>0$ and $c_8>0$ such that for all sufficiently large $n$,
$$
\P \Big\{ \sum_{y\in \G^{(L)}_n} {\bf 1}_{\{ \exists z\in \G^{(L)}_{2n+2}: \; y<z\} } \ge \ee^{c_7 \, n} \Big\}
\ge 
c_8.
$$  



\noindent Applying Claim \ref{claim:GW_embedded}, we see that there exists a constant $c_9>0$ such that for all sufficiently large $s$,
\begin{equation}
    \P \{ \# \mathscr{K}_s \ge \ee^{c_9 \, s} \}
    \ge 
    c_8.
    \label{sous-arbre-GW}
\end{equation}
 
Let $r>4s$. We have
$$
\gamma_r
:=
P_\omega ( T_{\mathscrHV_r} < T_{{\buildrel \leftarrow \over \varnothing}})
\ge
\sum_{x\in \mathscr{K}_s} P_\omega\{ T_{\mathscrHV_s} < T_{{\buildrel \leftarrow \over \varnothing}}, \; X_{T_{\mathscrHV_s}} =x\} \, \gamma_{r-V(x)}^{(x)} \, ,
$$

\noindent where, conditionally on $\mathscr{F}_{\!\! \mathscrHV_s}$, $(\gamma^{(x)}_t, \, t\ge 0)$, for $x\in \mathscr{K}_s$, are independent copies of $(\gamma_t, \, t\ge 0)$, and are independent of $\mathscr{F}_{\!\! \mathscrHV_s}$. [For $x\in \mathscr{K}_s$, we have $V(x) \le 4s <r$, so $\gamma_{r-V(x)}^{(x)}$ is well defined.] For $x\in \mathscr{K}_s$, and with the notation $\Lambda(x) := \sum_{y: \, {\buildrel \leftarrow \over y} = x} \ee^{-\Delta V(y)}$ from (\ref{Lambda}),
$$
P_\omega\{ T_{\mathscrHV_s} < T_{{\buildrel \leftarrow \over \varnothing}}, \; X_{T_{\mathscrHV_s}} =x\}
\ge
\prod_{j=1}^{|x|} \omega(x_{j-1}, \, x_j)
=
\frac{\ee^{-V(x)}}{\prod_{j=0}^{|x|-1} [1+\Lambda(x_j)]} \, ;
$$

\noindent on the other hand, by definition of $\mathscr{K}_s$, we have $\prod_{j=0}^{|x|-1} [1+\Lambda(x_j)]\le \ee^{2c_4L^{1-\alpha}s}$ and $V(x) \le 4s$ for $x\in \mathscr{K}_s$. Consequently, for $x\in \mathscr{K}_s$,
$$
P_\omega\{ T_{\mathscrHV_s} < T_{{\buildrel \leftarrow \over \varnothing}}, \; X_{T_{\mathscrHV_s}} =x\}
\ge 
\ee^{-(4+2c_4L^{1-\alpha})s} \, .
$$

\noindent Hence, writing $c_{10} := 4+2c_4L^{1-\alpha}$, we have
$$
\gamma_r
\ge
\ee^{-c_{10}\, s} \sum_{x\in \mathscr{K}_s} \gamma_{r-s}^{(x)} 
\ge
\ee^{-c_{10}\, s} \max_{x\in \mathscr{K}_s} \gamma_{r-s}^{(x)} .
$$

\noindent Applying (\ref{Paley-Zygmund-csq}) to $\gamma_{r-s}$ implies that if $r-s$ is sufficiently large,
\begin{eqnarray*}
    \P\{ \gamma_r \ge \ee^{-c_{10}\, s} \ee^{-(1+\varepsilon_2)(2(r-s))^{1/2}}\}
 &\ge&1- \E\{ (1-\ee^{-\varepsilon_2 (r-s)^{1/2}})^{\# \mathscr{K}_s} \}
    \\
 &\ge&1- \E\{ \ee^{-\ee^{-\varepsilon_2 (r-s)^{1/2}} \# \mathscr{K}_s} \}
    \\
 &\ge&(1-\ee^{-\ee^{-\varepsilon_2 (r-s)^{1/2}} \ee^{c_9 s}}) \, \P\{ \# \mathscr{K}_s \ge \ee^{c_9 s}\}.
\end{eqnarray*}

\noindent By (\ref{sous-arbre-GW}), $\P\{ \# \mathscr{K}_s \ge \ee^{c_9 s}\} \ge c_8$ if $s$ is sufficiently large. As a consequence, for all sufficiently large $s$ and $r-s$,
$$
\P\{ \gamma_r \ge \ee^{-c_{10}\, s} \ee^{-(1+\varepsilon_2)(2(r-s))^{1/2}}\}
\ge
c_8[1-\ee^{-\ee^{-\varepsilon_2 (r-s)^{1/2}} \ee^{c_9 s}}] .
$$ 

\noindent We take $s:= \frac{2}{c_9} \varepsilon_2\, r^{1/2}$, and see that for $\varepsilon_3 := (1+ \frac{2^{1/2} \, c_{10}}{c_9}) \varepsilon_2$, there exists $c_{11}\in (0, \, 1)$ such that for all sufficiently large $r$, say $r\ge r_0$,
\begin{equation}
    \P\{ \gamma_r \ge \ee^{-(1+\varepsilon_3)(2r)^{1/2}}\}
    \ge
    c_{11}.
    \label{0-1}
\end{equation}

Let $J_1$ be an integer such that $(1-c_{11})^{J_1}<\varepsilon_3$.
Let $\P^* (\, \cdot \, ) := \P ( \, \cdot \, | \, \hbox{\rm non-extinction})$ as before. Under $\P^*$, the system survives almost surely, so there exists an integer $J_2$ such that $\P^* \{ \sum_{|x|= J_2} 1 > J_1\} > 1-\varepsilon_3$. Let $r_1$ be sufficiently large such that $\P^* \{ \sum_{|x|= J_2} {\bf 1}_{\{ V(x) < r_1\} } \ge J_1\} \ge 1-\varepsilon_3$. We observe that for $r \ge r_1$,
\begin{eqnarray*}
    \gamma_r
 &\ge& \max_{y: \, |y| = J_2, \; V(y) < r_1} P_\omega\{ T_y < T_{{\buildrel \leftarrow \over \varnothing}} \} \, P_\omega^y \{ T_{\mathscrHV_r} < T_{{\buildrel \leftarrow \over \varnothing}} \}
    \\
 &\ge& c_{12}(\omega) \, \max_{y: \, |y| = J_2, \; V(y) < r_1} P_\omega^y \{ T_{\mathscrHV_r} < T_{{\buildrel \leftarrow \over \varnothing}} \} \, ,
\end{eqnarray*}

\noindent where $c_{12}(\omega) := \min_{y: \, |y| = J_2, \; V(y) < r_1} P_\omega\{ T_y < T_{{\buildrel \leftarrow \over \varnothing}} \} >0$ $\P$-a.s.\ (notation: $\min_{\varnothing} := 1$, $\max_{\varnothing} := 0$). 

For $|y| = J_2$ with $V(y) < r_1$, conditionally on $V(y)$, $P_\omega^y \{ T_{\mathscrHV_r} < T_{{\buildrel \leftarrow \over \varnothing}} \}$ is distributed as $\gamma_{r-V(y)}$, which is greater than or equal to $\gamma_r$. It follows from (\ref{0-1}) that for $r\ge \max\{ r_1, \, r_0\}$,
\begin{eqnarray*}
    \P\{ \gamma_r \ge c_{12}(\omega) \, \ee^{-(1+\varepsilon_3)(2r)^{1/2}}\}
 &\ge&\P\Big\{ \max_{y: \, |y| = J_2, \; V(y) < r_1} P_\omega^y \{ T_{\mathscrHV_r} < T_{{\buildrel \leftarrow \over \varnothing}} \}  \ge \ee^{-(1+\varepsilon_3)(2r)^{1/2}}\Big\}
    \\
 &\ge&(1-(1-c_{11})^{J_1}) \P \Big\{ \sum_{|x|= J_2} {\bf 1}_{\{ V(x) < r_1\} } \ge J_1 \Big\} \, .
\end{eqnarray*}

\noindent By definition of $r_1$, we have $\P \{ \sum_{|x|= J_2} {\bf 1}_{\{ V(x) < r_1\} } \ge J_1 \} \ge (1-\varepsilon_3) (1-q)$, where $q:= \P\{ \mathrm{extinction}\} < 1$. Therefore, for $r\ge \max\{ r_1, \, r_0\}$,
$$
\P\{ \gamma_r \ge c_{12}(\omega) \, \ee^{-(1+\varepsilon_3)(2r)^{1/2}}\}
\ge
(1-(1-c_{11})^{J_1}) (1-\varepsilon_3) (1-q)
\ge
(1-\varepsilon_3)^2 (1-q) \, ,
$$

\noindent the last inequality following from the definition of $J_1$. Since $c_{12}(\omega) >0$ $\P$-a.s., we have proved that
$$
\P^* \Big\{ \, \liminf_{r\to \infty} \, \frac{\log \gamma_r}{(2r)^{1/2}} \ge -1-\varepsilon_3 \, \Big\} 
\ge 
(1-\varepsilon_3)^2\,  .
$$

\noindent Recall the definition $\varepsilon_3 := (1+ \frac{2^{1/2} \, c_{10}}{c_9}) \varepsilon_2$, with $\varepsilon_2  > 6\varepsilon + 8\varepsilon_1$, $\varepsilon>0$ and $\varepsilon_1 \in (0, \, c_1 \, \varepsilon)$; so $\varepsilon_3>0$ can be taken arbitrarily small. This yields the lower bound in (\ref{key-estimate}), and thus completes the proof of Theorem \ref{t:main} by admitting Lemmas \ref{l:Sigma} and \ref{l:Z}.\hfill$\Box$

\medskip

The rest of the section is devoted to the proof of Lemmas \ref{l:Sigma} and \ref{l:Z}.

\subsection{Proof of Lemma \ref{l:Sigma}}
\label{subs:proof-first-lemma}

In the study of one-dimensional random walks, a frequent type of technical difficulties is to handle the overshoots. Such difficulties are, unfortunately, present throughout the proof of both Lemmas \ref{l:Sigma} and \ref{l:Z}.    

Let $r>0$. Let $\chi \in (0, \, 1)$. Recall from (\ref{k})--(\ref{hm}) that
$$
    k := \lfloor r^{1-\chi} \rfloor \, ,
    \qquad
    h_m
    :=
    \frac{r}{k} \, m \, ,
    \qquad 
    0\le m\le k \, .
$$

\noindent Recall from (\ref{M}) that $\mathscrHV_r := \{ x\in \T: \, V(x) \ge r, \; \overline{V}({\buildrel \leftarrow \over x}) < r\}$. We distinguish the vertices $x$ of $\mathscrHV_r$ according to whether there are some ``large overshoots" of the random potential $V(\cdot)$ along the path $[\![ \varnothing, \, x]\!]$: let $\theta \in (\frac12, \, \chi)$, and let
\begin{eqnarray*}
    \mathscrHV_{r, \, +} 
 &:=& \Big\{ x\in \mathscrHV_r: \, \max_{1\le m< k} \Delta V(x_{\mathscrHx_{h_m}}) > r^\theta \Big\} ,
    \\
    \mathscrHV_{r, \, -} 
 &:=&\Big\{ x\in \mathscrHV_r: \, \max_{1\le m< k} \Delta V(x_{\mathscrHx_{h_m}}) \le r^\theta \Big\} ,
\end{eqnarray*}

\noindent where, as before, $\Delta V(y) := V(y) - V({\buildrel \leftarrow \over y})$ for any vertex $y\in \T \backslash \{ \varnothing\}$. 

Recall from (\ref{TMr}) that 
$$
T_{\mathscrHV_r} 
= 
\inf_{x\in \mathscrHV_r} T_x
=
\min\Big\{ \inf_{x\in \mathscrHV_{r, \, +} } T_x, \, \inf_{x\in \mathscrHV_{r, \, -} } T_x \Big\} \, ,
$$

\noindent where $T_x := \inf \{ i\ge 0: \, X_i =x\}$ as in (\ref{Tx}). So
\begin{equation}
    P_\omega ( T_{\mathscrHV_r} < T_{{\buildrel \leftarrow \over \varnothing}})
    \le
    \sum_{x\in \mathscrHV_{r, \, +} } P_\omega( T_x < T_{{\buildrel \leftarrow \over \varnothing}})
    +
    P_\omega \Big( \inf_{x\in \mathscrHV_{r, \, -} } T_x < T_{{\buildrel \leftarrow \over \varnothing}} \Big) \, .
    \label{proof-lemma4.1-start}
\end{equation}

We first bound $\sum_{x\in \mathscrHV_{r, \, +} } P_\omega( T_x < T_{{\buildrel \leftarrow \over \varnothing}})$.  By a one-dimensional argument (Golosov~\cite{golosov}), for any $x$, $y\in \T$ with $y<x$, 
\begin{equation}
    P_\omega\{ T_x < T_{{\buildrel \leftarrow \over \varnothing}} \, | \, X_0=y\} 
    = 
    \frac{\sum_{u\in [\![ \varnothing, \, y]\!]} \ee^{V(u)}}{\sum_{u\in [\![ \varnothing, \, x]\!]} \ee^{V(u)}} .
    \label{P(T_x<T_0):cas_general}
\end{equation}

\noindent In particular, for any $x\in \T \backslash \{ \varnothing\}$,
\begin{equation}
    P_\omega\{ T_x < T_{{\buildrel \leftarrow \over \varnothing}}\} 
    = 
    \frac{1}{\sum_{u\in [\![ \varnothing, \, x]\!]} \ee^{V(u)}}
    \le
    \ee^{- \overline{V}(x)} .
    \label{P(T_x<T_0)=}
\end{equation}

\noindent Hence
$$
\sum_{x\in \mathscrHV_{r, \, +} } P_\omega( T_x < T_{{\buildrel \leftarrow \over \varnothing}})
\le
\sum_{x\in \mathscrHV_{r, \, +}} \ee^{- \overline{V}(x)}
=
\sum_{x\in \mathscrHV_{r, \, +}} \ee^{- V(x)} ,
$$

\noindent the last identity following from the fact that $\overline{V}(x) = V(x)$ for all $x\in \mathscrHV_{r, \, +}$. Taking expectation with respect to $\E$ on both sides, we obtain, by means of (\ref{many-to-one:stopping-lines}),
$$
\E\Big[ \sum_{x\in \mathscrHV_{r, \, +} } P_\omega( T_x < T_{{\buildrel \leftarrow \over \varnothing}}) \Big]
\le
\Q \Big[ \max_{1\le m < k} \Delta S_{\mathscrHS_{h_m}} > r^\theta \Big]
\le
\sum_{m=1}^{k-1} \Q \Big[ \Delta S_{\mathscrHS_{h_m}} > r^\theta \Big] .
$$

\noindent We use (\ref{moment-expo-overshoot}) to see that for some constant $c_{13}>0$,
$$
    \E\Big[ \sum_{x\in \mathscrHV_{r, \, +} } P_\omega( T_x < T_{{\buildrel \leftarrow \over \varnothing}}) \Big]
    \le
    c_{13} \, (k-1) \, \ee^{-c_2 \, r^\theta}
    =
    c_{13} \, (\lfloor r^{1-\chi} \rfloor-1) \, \ee^{-c_2 \, r^\theta} .
$$

\noindent Recall that $\theta> \frac12$. In view of (\ref{proof-lemma4.1-start}), the proof of Lemma \ref{l:Sigma} is reduced to showing the following:
\begin{equation}
    \limsup_{r\to \infty} \frac{1}{(2r)^{1/2}} \, \log \E \Big[ P_\omega \Big( \inf_{x\in \mathscrHV_{r, \, -} } T_x < T_{{\buildrel \leftarrow \over \varnothing}} \Big) \Big]
    \le
    -1 \, .
    \label{proof-lemma4.1-step1}
\end{equation}

For any vertex $x\in \mathscrHV_r$, let us recall $a_j^{(x)}$ from (\ref{aV}), and define
$$
    \tau_x :=
    \inf \{ j: \, 1 \le j \le |x|, \, \overline{V}(x_j) - V(x_j) \ge a_j^{(x)} \}.
    \qquad (\inf\varnothing := \infty)
$$

\noindent For $x\in \mathscrHV_r$, we have either $\tau_x < |x|$ (with strict inequality), or $\tau_x = \infty$. We observe that
\begin{eqnarray*}
    \inf_{x\in \mathscrHV_{r, \, -} } T_x 
 &=&\min\Big\{ \inf_{x\in \mathscrHV_{r, \, -}: \, \tau_x < |x|} T_x, \; \inf_{x\in \mathscrHV_{r, \, -}: \, \tau_x =\infty} T_x \Big\}  
    \\
 &\ge&\min\Big\{ \inf_{x\in \mathscrHV_{r, \, -}: \, \tau_x < |x|} T_{y(x)}, \; \inf_{x\in \mathscrHV_{r, \, -}: \, \tau_x =\infty} T_x \Big\} ,
\end{eqnarray*}

\noindent where $y(x) := x_{\tau_x}$. Hence
\begin{eqnarray}
    P_\omega \Big( \inf_{x\in \mathscrHV_{r, \, -} } T_x < T_{{\buildrel \leftarrow \over \varnothing}} \Big)
 &\le&P_\omega\Big( \inf_{x\in \mathscrHV_{r, \, -}: \, \tau_x < |x|} T_{y(x)} < T_{{\buildrel \leftarrow \over \varnothing}} \Big)
  + \sum_{x\in \mathscrHV_{r, \, -}: \, \tau_x =\infty} P_\omega\{ T_x < T_{{\buildrel \leftarrow \over \varnothing}}\} 
    \nonumber
    \\
 &=:& \Sigma_1 + \Sigma_2 \, ,
    \label{Sigma1-2}
\end{eqnarray}

\noindent with obvious notation. It is easy to get an upper bound for $\Sigma_2$: by (\ref{P(T_x<T_0)=}), $P_\omega\{ T_x < T_{{\buildrel \leftarrow \over \varnothing}}\} \le \ee^{- \overline{V}(x)}$ (which is $\ee^{-V(x)}$ for $x\in \mathscrHV_{r, \, -}$), whereas $\tau_x =\infty$ implies $\overline{V}(x_i) - V(x_i) < a_i^{(x)}, \; \forall i< |x|$, so
\begin{equation}
    \Sigma_2
    \le
    \sum_{x\in \mathscrHV_r} \ee^{-V(x)} \, {\bf 1}_{\{ \overline{V}(x_i) - V(x_i) < a_i^{(x)}, \; \forall i< |x|\} } \prod_{m=1}^{k-1} {\bf 1}_{\{ \Delta V(x_{\mathscrHx_{h_m}}) \le r^\theta \} } \, .
    \label{Sigma_2}
\end{equation}

\noindent To bound $\Sigma_1$, we note that
$$
\inf_{x\in \mathscrHV_{r, \, -}: \, \tau_x < |x|} T_{y(x)}
=
\inf\{ T_y: \; \exists x\in \mathscrHV_{r, \, -}, \; y=x_{\tau_x}, \; \tau_x < |x| \}  \, .
$$

\noindent Let $y\in \T$ with $j:= |y| \ge 1$ such that $h_{m-1} \le \overline{V}(y) < h_m$ for some $m \in [1, \, k]$. We define
$$
a_i^{(y)}
:=
\begin{cases}
    \lambda_\ell \, , &  \hbox{\rm if $\mathscrHy_{h_{\ell-1}} \le i < \mathscrHy_{h_\ell} \, $ for $\ell \in [1, \, m)$} \, ,
    \\
    \lambda_m \, , & \hbox{\rm if $\mathscrHy_{h_{m-1}} \le i \le j$} \, .
\end{cases}
$$

\noindent Clearly, if $y=x_{\tau_x}$ for some $x\in \mathscrHV_{r, \, -}$ satisfying $\tau_x < |x|$, then $\overline{V}(y_i) - V(y_i) < a_i^{(y)}$, $\forall i<j$, and $\overline{V}(y_j) - V(y_j) \ge \lambda_m$, and moreover $\Delta V(y_{\mathscrHy_{h_\ell}}) \le r^\theta$, $\forall 1\le \ell < m$. Accordingly,
\begin{eqnarray*}
    \Sigma_1
 &\le&\sum_{m=1}^k \sum_{j=1}^\infty \sum_{|y|=j}
    {\bf 1}_{\{ h_{m-1} \le \overline{V}(y) < h_m \} } \, 
    {\bf 1}_{\{ \overline{V}(y_i) - V(y_i) < a_i^{(y)}, \; \forall i<j; \; \overline{V}(y_j) - V(y_j) \ge \lambda_m\} } \times
    \\
 &&\qquad \times \Big( \prod_{\ell=1}^{m-1} {\bf 1}_{\{ \Delta V(y_{\mathscrHy_{h_\ell}}) \le r^\theta \} } \Big) \,
    P_\omega\{ T_y < T_{{\buildrel \leftarrow \over \varnothing}}\} \, .
\end{eqnarray*}

\noindent Again, by (\ref{P(T_x<T_0)=}), we have $P_\omega\{ T_y < T_{{\buildrel \leftarrow \over \varnothing}}\} \le \ee^{-\overline{V}(y)}$. This gives the analogue of (\ref{Sigma_2}) for $\Sigma_1$. 

We apply the many-to-one formula in (\ref{many-to-one}). Recall from (\ref{mathscrHS_r}) that $\mathscrHS_u := \inf\{ i\ge 0: \, S_i \ge u \}$ (for $u\ge 0$), and from (\ref{DeltaS}) that $\Delta S_i := S_i -S_{i-1}$. Define
\begin{equation}
    a_i^{(S)} 
    := 
    \lambda_m,
    \qquad
    \hbox{\rm if $\mathscrHS_{h_{m-1}} \le i< \mathscrHS_{h_m}$ and $1\le m\le k$} \, .
    \label{aS}
\end{equation}

\noindent By (\ref{many-to-one}), 
\begin{eqnarray}
    \E(\Sigma_1)
 &\le&\sum_{m=1}^k \sum_{j=1}^\infty 
    \E_{\Q}\Big[ \ee^{-(\overline{S}_j -S_j)} \, 
    {\bf 1}_{\{ h_{m-1} \le \overline{S}_j < h_m \} } \, 
    {\bf 1}_{\{ \overline{S}_i - S_i < a_i^{(S)}, \; \forall i<j; \; \overline{S}_j - S_j \ge \lambda_m\} }\times
    \nonumber
    \\
 &&\qquad\qquad \times
    \prod_{\ell=1}^{m-1} {\bf 1}_{\{ \Delta S_{\mathscrHS_{h_\ell}} \le r^\theta \} } \Big]
    \nonumber
    \\
 &\le&\sum_{m=1}^k \sum_{j=1}^\infty \ee^{-\lambda_m}\, 
    \E_{\Q}\Big[ {\bf 1}_{\{ h_{m-1} \le \overline{S}_j < h_m \} } \, 
    {\bf 1}_{\{ \overline{S}_i - S_i < a_i^{(S)}, \; \forall i<j\} }\,
    \prod_{\ell=1}^{m-1} {\bf 1}_{\{ \Delta S_{\mathscrHS_{h_\ell}} \le r^\theta \} } \Big] \, .
    \label{E(Sigma_1)<}
\end{eqnarray}

\noindent Similarly, applying (\ref{many-to-one:stopping-lines}) in place of (\ref{many-to-one}) to $\E(\Sigma_2)$, we obtain:
\begin{equation}
    \E(\Sigma_2)
    \le
    \Q \Big\{ \overline{S}_i - S_i < a_i^{(S)}, \; \forall 1\le i < \mathscrHS_r; \; \max_{1\le \ell <k} \Delta S_{\mathscrHS_{h_\ell}} \le r^\theta\Big\}.
    \label{E(Sigma_2)<}
\end{equation}

\noindent At this stage, we have two preliminary results.

\begin{claim}
\label{claim:proba}

 For any integers $1\le m_0 \le m <k$ and any $s\in (-\infty, \, h_{m_0})$, we
 define
 \begin{equation}
     f_{m_0,m}(s)
     :=
     \Q \Big( \bigcap_{\ell=m_0+1}^{m+1} \{ \max_{i\in [\mathscrHS_{h_{\ell-1}-s}, \, \mathscrHS_{h_\ell-s})}(\overline{S}_i - S_i) < \lambda_\ell\} \cap \bigcap_{\ell=m_0}^m \{ \Delta S_{\mathscrHS_{h_\ell-s}} \le r^\theta \}  \Big) .
     \label{f}
 \end{equation}
 Then, as $r\to \infty$,
 \begin{equation}
     \sup_{s<h_{m_0}}f_{m_0,m}(s)
     \le
     \ee^{- (1+o(1)) \sum_{\ell=m_0+1}^{m+1} \frac{r^\chi}{\lambda_\ell}} \, ,
     \label{f<}
 \end{equation}
 uniformly in $1\le m_0 \le m <k$. Furthermore,
 \begin{equation}
     \Q \Big( \bigcap_{\ell=1}^{m+1} \{ \max_{i\in [\mathscrHS_{h_{\ell-1}}, \, \mathscrHS_{h_\ell})}(\overline{S}_i - S_i) < \lambda_\ell\} \cap \bigcap_{\ell=1}^m \{ \Delta S_{\mathscrHS_{h_\ell}} \le r^\theta \}  \Big)
     \le
     \ee^{- (1+o(1)) \sum_{\ell=1}^{m+1} \frac{r^\chi}{\lambda_\ell}} \, ,
     \label{f<:m=0}
 \end{equation}
 uniformly in $1\le m<k$. 
 
\end{claim}

\begin{claim}
\label{claim:sum_j}

 There exists a constant $c_{14}>0$ such that for $r\to \infty$,
 \begin{eqnarray}
  &&\sum_{j=1}^\infty 
     \E_{\Q}\Big[ {\bf 1}_{\{ h_{m-1} \le \overline{S}_j < h_m \} } \, 
     {\bf 1}_{\{ \overline{S}_i - S_i < a_i^{(S)}, \; \forall i<j\} }\,
     \prod_{\ell=1}^{m-1} {\bf 1}_{\{ \Delta S_{\mathscrHS_{h_\ell}} \le r^\theta \} } \Big]
     \nonumber
     \\
  &\le&c_{14} \, r \, \exp \Big( - (1+o(1))\sum_{\ell=1}^{m-1} \frac{r^\chi}{\lambda_\ell} \Big) ,
     \label{sum_j}
 \end{eqnarray}
 uniformly in $m\in [1, \, k]$.

\end{claim}

\medskip

\noindent {\it Proof of Claim \ref{claim:proba}.} Applying the strong Markov property successively at $\mathscrHS_{h_m-s}$, $\mathscrHS_{h_{m-1}-s}$, $\cdots$, $\mathscrHS_{h_{m_0}-s}$, we obtain that
$$
f_{m_0,m}(s)
\le
\prod_{\ell=m_0+1}^{m+1} \sup_{u\in [0, \, r^\theta]} \Q\Big( \max_{0\le i< \mathscrHS_{h_\ell- h_{\ell-1}-u} } (\overline{S}_i -S_i) < \lambda_\ell \Big).
$$

\noindent By Lemma \ref{l:RW-chute-ub}, we arrive at the following estimate: when $r\to \infty$,
$$
f_{m_0,m}(s)
\le
\exp\Big( - (1+o(1)) \sum_{\ell=m_0+1}^{m+1} \frac{h_\ell- h_{\ell-1}-r^\theta}{\lambda_\ell} \Big)
\le
\exp\Big( - (1+o(1)) \sum_{\ell=m_0+1}^{m+1} \frac{r^\chi}{\lambda_\ell} \Big) ,
$$

\noindent uniformly in $s<h_{m_0}$ and in $1\le m_0\le m <k$;\footnote{Since $h_m-h_{m-1} = \frac{r}{k}$ (by (\ref{hm})), it is here we use the condition $\theta < \chi$ to ensure $h_m-h_{m-1}-r^\theta>0$.} this yields (\ref{f<}). The proof of (\ref{f<:m=0}) is along the same lines.\hfill$\Box$

\bigskip

\noindent {\it Proof of Claim \ref{claim:sum_j}.} Let $\mathrm{LHS}_{(\ref{sum_j})}$ denote the sum on the left-hand side of (\ref{sum_j}). Then
$$
\mathrm{LHS}_{(\ref{sum_j})}
=
\E_\Q \Big[ 
    \Big( \prod_{\ell=1}^{m-1} {\bf 1}_{\{ \Delta S_{\mathscrHS_{h_\ell}} \le r^\theta \} } \Big)
    \sum_{j=\mathscrHS_{h_{m-1}}}^{\mathscrHS_{h_m} - 1} 
    {\bf 1}_{\{ \overline{S}_i - S_i < a_i^{(S)}, \; \forall i<j\} } 
    \Big] .
$$

\noindent By definition of $a_i^{(S)}$ in (\ref{aS}), this yields that
\begin{eqnarray*}
    \mathrm{LHS}_{(\ref{sum_j})}
 &=& \E_\Q \Big[ 
    \sum_{j=\mathscrHS_{h_{m-1}}}^{\mathscrHS_{h_m} - 1} 
    {\bf 1}_{\{ \overline{S}_i - S_i < \lambda_m, \, \forall i\in [\mathscrHS_{h_{m-1}}, \, j)\} } \times
    \\
 &&\qquad \times
    \Big( \prod_{\ell=1}^{m-1} {\bf 1}_{\{ \max_{i\in [\mathscrHS_{h_{\ell-1}}, \, \mathscrHS_{h_\ell})}(\overline{S}_i - S_i) < \lambda_\ell\} \cap \{ \Delta S_{\mathscrHS_{h_\ell}} \le r^\theta \} } \Big) \Big] .
\end{eqnarray*}

\noindent We proceed to get rid of the sum over $j$ on the right-hand side. Applying the strong Markov property at time $\mathscrHS_{h_{m-1}}$, we have
\begin{equation}
    \mathrm{LHS}_{(\ref{sum_j})}
    \le
    \E_\Q \Big[ \Big( \prod_{\ell=1}^{m-1} {\bf 1}_{\{ \max_{i\in [\mathscrHS_{h_{\ell-1}}, \, \mathscrHS_{h_\ell})}(\overline{S}_i - S_i) < \lambda_\ell\} \cap \{ \Delta S_{\mathscrHS_{h_\ell}} \le r^\theta \} } \Big) \, \Xi_m  \Big],
    \label{proof_claim_sum_j-1}
\end{equation}

\noindent where
\begin{eqnarray*}
    \Xi_m
 &:=&\sup_{x\in [h_m-h_{m-1}-r^\theta, \, h_m-h_{m-1}]} \E_\Q\Big( \sum_{j=0}^{\mathscrHS_x - 1} {\bf 1}_{\{ \overline{S}_i - S_i < \lambda_m, \, \forall i\in [0, \, j)\} } \Big)
    \\
 &\le&\E_\Q\Big( \sum_{j=0}^\infty {\bf 1}_{\{ \overline{S}_i - S_i < \lambda_m, \, \forall i\in [0, \, j) \} } \Big)  \, .
\end{eqnarray*}

\noindent To estimate the expectation on the right-hand side, we write $\sum_{j=0}^\infty = \sum_{n=1}^\infty \sum_{j=(n-1) \lambda_m^2}^{n \lambda_m^2 -1}$ (by implicitly treating $\lambda_m^2$ as an integer; otherwise we replace $\lambda_m$ by $\lceil \lambda_m \rceil$, and the next three paragraphs will still go through with obvious modifications), so that
\begin{eqnarray*}
    \Xi_m
 &\le&\sum_{n=1}^\infty \E_\Q\Big( \sum_{j=(n-1) \lambda_m^2}^{n \lambda_m^2 -1} {\bf 1}_{\{ \overline{S}_i - S_i < \lambda_m, \, \forall i\in [0, \, j) \} } \Big)
    \\
 &\le&\sum_{n=1}^\infty \lambda_m^2\, \Q \Big\{ \max_{0\le i< (n-1) \lambda_m^2} ( \overline{S}_i - S_i) < \lambda_m \Big\} \, .
\end{eqnarray*}

\noindent By the Markov property, $\Q \{ \max_{0\le i< (n-1) \lambda_m^2} ( \overline{S}_i - S_i) < \lambda_m \} \le [ \, \Q \{ \max_{0\le i< \lambda_m^2} ( \overline{S}_i - S_i) < \lambda_m \}]^{n-1}$. So
$$
\Xi_m
\le
\sum_{n=1}^\infty \lambda_m^2 \Big[ \, \Q \Big\{ \max_{0\le i< \lambda_m^2} ( \overline{S}_i - S_i) < \lambda_m \Big\} \Big]^{n-1} \, .
$$

\noindent We let $r\to \infty$ (so that $\lambda_m \to \infty$ uniformly in $m\in [1, \, k]$). By Donsker's theorem, $\Q \{ \max_{0\le i< \lambda_m^2} ( \overline{S}_i - S_i) < \lambda_m \} \to \p \{ \sup_{s\in [0, \, 1]} (\overline{W}_s-W_s) <\frac{1}{\sigma}\} < 1$, where $(W_s, \, s\ge 0)$ under $\p$ is a standard Brownian motion, and $\overline{W}_s := \sup_{u\in [0, \, s]} W_u$. So there exists a constant $0<c_{15}<1$ such that for all sufficiently large $r$ and all $m\in [1, \, k]$, $\Q \{ \max_{0\le i< \lambda_m^2} ( \overline{S}_i - S_i) < \lambda_m \} \le 1-c_{15}$, which, in turn, yields that 
$$
\Xi_m
\le
\sum_{n=1}^\infty \lambda_m^2 (1-c_{15})^{n-1}
=
\frac{\lambda_m^2}{c_{15}}
\le
\frac{2r}{c_{15}} \, .
$$

\noindent Going back to (\ref{proof_claim_sum_j-1}), this yields that for all sufficiently large $r$ (writing $c_{16} := \frac{2}{c_{15}}$),
$$
    \mathrm{LHS}_{(\ref{sum_j})}
    \le
    c_{16} \, r \, \Q \Big( \bigcap_{\ell=1}^{m-1} \{ \max_{i\in [\mathscrHS_{h_{\ell-1}}, \, \mathscrHS_{h_\ell})}(\overline{S}_i - S_i) < \lambda_\ell\} \cap \{ \Delta S_{\mathscrHS_{h_\ell}} \le r^\theta \} \Big) .
$$

\noindent This implies Claim \ref{claim:sum_j} in case $2\le m<k$ by means of (\ref{f<:m=0}), and trivially in case $m=1$.\hfill$\Box$

\bigskip

We continue with the proof of Lemma \ref{l:Sigma}. By (\ref{E(Sigma_1)<}) and Claim \ref{claim:sum_j}, we have
$$
\E(\Sigma_1)
\le
c_{14} \, r \, \sum_{m=1}^k \exp \Big( - \lambda_m - (1+o(1)) \sum_{\ell=1}^{m-1} \frac{r^\chi}{\lambda_\ell} \Big) .
$$

\noindent By definition, $k := \lfloor r^{1-\chi}\rfloor$ and $\lambda_m:= (2r)^{1/2} \, (\frac{k-m+1}{k})^{1/2}$; hence for $r\to \infty$,
\begin{equation}
    \sum_{\ell=1}^{m+1} \frac{r^\chi}{\lambda_\ell}
    =
    (2r^\chi)^{1/2}[k^{1/2} - (k-m)^{1/2}] 
    +
    o((2r)^{1/2})\, ,
    \label{sum_over_lambda_j}
\end{equation}

\noindent uniformly in $1\le m_0\le m <k$. In particular,
\begin{equation}
    \sum_{\ell=1}^k \frac{r^\chi}{\lambda_\ell}
    \sim
    (2r)^{1/2}\, .
    \label{sum(hm-hm-1)/lambda_m=sqrt(2r)}
\end{equation}

\noindent So uniformly in $m\in [1, \,k]$,
\begin{eqnarray*}
&&\lambda_m + (2r^\chi)^{1/2}[k^{1/2} - (k-m+1)^{1/2}]
    \\
 &\ge&(1+o(1)) (2r)^{1/2} \inf_{s\in [0, \, 1]} \Big( (1-s)^{1/2} + [1 - (1-s)^{1/2}] \Big) \, ,
\end{eqnarray*}

\noindent and the infimum equals $1$ because the function $s\mapsto (1-s)^{1/2} + [1 - (1-s)^{1/2}]$ is identically $1$ on $[0, \, 1]$. Therefore,
$$
\E(\Sigma_1)
\le
c_{14} \, r k \ee^{-(1+o(1)) (2r)^{1/2}}
\le
\ee^{-(1+o(1)) (2r)^{1/2}} \, ,
$$

\noindent the second inequality being a consequence of definition $k := \lfloor r^{1-\chi}\rfloor$.

On the other hand, by (\ref{E(Sigma_2)<}) and (\ref{f<:m=0}) (applied to $m:=k-1$), we have 
$$
\E(\Sigma_2)
\le
\ee^{- (1+o(1)) \sum_{\ell=1}^k \frac{r^\chi}{\lambda_\ell}}
\le
\ee^{- (1+o(1)) (2r)^{1/2}} \, ,
$$

\noindent the second inequality being a consequence of (\ref{sum_over_lambda_j}) (applied to $m:= k-1$). Since $P_\omega ( \inf_{x\in \mathscrHV_{r, \, -} } T_x < T_{{\buildrel \leftarrow \over \varnothing}} ) \le \Sigma_1 + \Sigma_2$ (see (\ref{Sigma1-2})), this yields (\ref{proof-lemma4.1-step1}), and completes the proof of Lemma \ref{l:Sigma}.\hfill$\Box$ 

\medskip

The rest of the section is devoted to the proof of Lemma \ref{l:Z}, which is more technical. For the sake of clarity, we prove the three parts (\ref{E(Z)}), (\ref{E(Z2)}) and (\ref{(EZ)2}) separately.

\subsection{Proof of Lemma \ref{l:Z}: the inequality (\ref{E(Z)})}
\label{subs:proof-second-lemma-(i)}

Recall from (\ref{Z}) the definition $Z_r := \sum_{x\in \mathscrHV_r^*} {\bf 1}_{\{ T_x < T_{{\buildrel \leftarrow \over \varnothing}} \} }$, where 
\begin{eqnarray*}
    \mathscrHV_r^*
 &:=& \Big\{ x\in \mathscrHV_r: \, \max_{1\le m<k} \Delta V(x_{\mathscrHx_{h_m}}) \le r^\theta , \; \underline{V}(x) \ge - \beta , \; |x| < \lfloor \ee^{\varepsilon_1 \, r^{1/2}}\rfloor \, , \;
    \\
 && \overline{V}(x_j) - V(x_j) \le a_j^{(x)} , \; \forall 0\le j<|x| , \;  \max_{0\le \ell<|x|} \Lambda(x_\ell) \le \ee^{\varepsilon r^{1/2}} \Big\} ,
\end{eqnarray*}

\noindent with $\Lambda(x) := \sum_{y: \, {\buildrel \leftarrow \over y} = x} \ee^{-\Delta V(y)}$ as in (\ref{Lambda}). For brevity, we write, in this subsection,
$$
n
=
n(\varepsilon_1, \, r)
:=
\lfloor \ee^{\varepsilon_1\, r^{1/2}} \rfloor \, ;
$$

\noindent so $|x|+1 \le n$ for all $x\in \mathscrHV_r^*$. Since only $T_x$ and $T_{{\buildrel \leftarrow \over \varnothing}}$ depend on the biased walk $(X_i)$, we have
\begin{equation}
    E_\omega(Z_r)
    =
    \sum_{x\in \mathscrHV_r^*} 
    P_\omega \{ T_x < T_{{\buildrel \leftarrow \over \varnothing}} \} \, .
    \label{E(Z)=}
\end{equation}

\noindent By the identity in (\ref{P(T_x<T_0)=}), we have $P_\omega\{ T_x < T_{{\buildrel \leftarrow \over \varnothing}}\} \ge \frac{1}{|x|+1} \ee^{- \overline{V}(x)}$, which is $\ge \frac1n \ee^{- \overline{V}(x)} = \frac1n \ee^{- V (x)} $ for all $x\in \mathscrHV_r^*$. Taking expectation with respect to $\E$ on both sides leads to:
\begin{eqnarray*}
    \E[ E_\omega(Z_r) 
    ]
 &\ge&\frac1n \, \E \Big[ 
    \sum_{x\in \mathscrHV_r^*} \ee^{- V(x)} \Big] 
    \\
 &=&\frac1n \, \E \Big[ 
    \sum_{x\in \mathscrHV_r} 
    \ee^{- V(x)} \, 
    {\bf 1}_{\{ \overline{V}(x_j) - V(x_j) < a^{(x)}_j, \, \forall 0\le j < |x| \} } \, 
    {\bf 1}_{\{ \underline{V}(x) \ge -\beta\} } \times
    \\
 && \qquad\times {\bf 1}_{\{ |x|<n\} } \,
    {\bf 1}_{\{ \Lambda(x_\ell) \le \ee^{\varepsilon r^{1/2}}, \, \forall 0\le \ell< |x|\} } 
    \prod_{m=1}^{k-1} {\bf 1}_{\{ \Delta V(x_{\mathscrHx_{h_m}}) \le r^\theta\} } \Big] .
\end{eqnarray*}

\noindent The expression on the right-hand side is, according to formula (\ref{decomposition-epinale:stopping-lines}),
\begin{eqnarray*}
 &=& \frac{1}{n} \, 
    \Q \Big[ 
    \bigcap_{j=0}^{\mathscrHS_r-1} \{ \overline{S}_j - S_j < a^{(S)}_j, \, 
    S_j \ge - \beta \} 
    \cap 
    \nonumber 
    \\
 && \qquad \cap
    \{ \mathscrHS_r < n\} 
    \cap 
    \bigcap_{\ell=1}^{\mathscrHS_r} 
    \{ \eta_\ell \le \ee^{\varepsilon r^{1/2}}\} 
    \cap 
    \bigcap_{m=1}^{k-1} 
    \{ \Delta S_{\mathscrHS_{h_m}} \le r^\theta \} \Big] ,
\end{eqnarray*}

\noindent where $\mathscrHS_r := \inf\{ i\ge 0: \, S_i \ge r \}$ as in (\ref{mathscrHS_r}), $\overline{S}_j := \max_{0\le i\le j} S_i$ as in (\ref{S_bar}), $\Delta S_j := S_j-S_{j-1}$ as before (with $S_0 := 0$), and $\eta_\ell := \sum_{y: \, {\buildrel \leftarrow \over y} = w_{\ell-1}} \ee^{-\Delta V(y)}$. [In particular, $\eta_1 := \sum_{y: \, |y|=1} \ee^{- V(y)}$.] Recall from Section \ref{s:preliminaries} that $(\Delta S_i, \, \eta_i)$, $i\ge 1$, are i.i.d.\ random vectors under $\Q$. Hence
\begin{equation}
    \E[ E_\omega(Z_r) ]
    \ge
    \frac{1}{n} \, [q_1(r) - q_2(r)]\, ,
    \label{E(Z)>q1-q2}
\end{equation}

\noindent where
\begin{eqnarray*}
    q_1(r)
 &:=& \Q \Big[ 
    \bigcap_{j=0}^{\mathscrHS_r-1} \{ \overline{S}_j - S_j < a^{(S)}_j, \, 
    S_j \ge - \beta \} 
    \cap 
    \{ \mathscrHS_r < n \} 
    \cap 
    \bigcap_{m=1}^{k-1}
    \{ \Delta S_{\mathscrHS_{h_m}} \le r^\theta \} 
    \Big] ,
    \\
    q_2(r)
 &:=& \Q \Big[ 
    \bigcap_{j=0}^{\mathscrHS_r-1} \{ \overline{S}_j - S_j < a^{(S)}_j\}     
    \cap 
    \bigcap_{m=1}^{k-1}
    \{ \Delta S_{\mathscrHS_{h_m}} \le r^\theta \} 
    \cap 
    \bigcup_{\ell=1}^{\mathscrHS_r \wedge n} \{ \eta_\ell > \ee^{\varepsilon r^{1/2}} \} 
    \Big] .
\end{eqnarray*}


\noindent By definition of $(a^{(S)}_j)$ in (\ref{aS}) (with notation $\Delta S_0:=0$ for the term $m=1$ below),
\begin{eqnarray*}
    q_1(r)
 &=& \Q \Big( \{ \mathscrHS_r < n\} 
    \cap
    \\
 &&\cap \bigcap_{m=1}^k \bigcap_{j=\mathscrHS_{h_{m-1}}}^{\mathscrHS_{h_m}-1}
    \{ \overline{S}_j - S_j < \lambda_m, \, 
    S_j \ge - \beta\} \cap \{ \Delta S_{\mathscrHS_{h_{m-1}}} \le r^\theta \} \Big) .
\end{eqnarray*}

\noindent Since $\{ \mathscrHS_r < n \} \supset \cap_{m=1}^k \{ \mathscrHS_{h_m} - \mathscrHS_{h_{m-1}} < \lfloor \frac{n}{k} \rfloor\}$, we have
\begin{eqnarray*}
    q_1(r)
 &\ge& \Q \Big\{ 
    \bigcap_{m=1}^k \bigcap_{j=\mathscrHS_{h_{m-1}}}^{\mathscrHS_{h_m}-1}
    \{ 
    \overline{S}_j - S_j < \lambda_m, \, 
    S_j \ge - \beta\} \cap
    \\
 &&\qquad\qquad 
    \cap \{ \Delta S_{\mathscrHS_{h_{m-1}}} \le r^\theta, \, 
    \mathscrHS_{h_m} - \mathscrHS_{h_{m-1}} < \lfloor \frac{n}{k} \rfloor \} 
    \Big\}
    \\
 &\ge& \Q \Big\{ 
    \bigcap_{m=1}^k \bigcap_{j=\mathscrHS_{h_{m-1}}}^{\mathscrHS_{h_m}-1}
    \{ 
    \overline{S}_j - S_j < \lambda_m, \, 
    S_j - S_{\mathscrHS_{h_{m-1}}} \ge - \beta \} \cap
    \\
 &&\qquad\qquad 
    \cap \{ \Delta S_{\mathscrHS_{h_{m-1}}} \le r^\theta, \, 
    \mathscrHS_{h_m} - \mathscrHS_{h_{m-1}} < \lfloor \frac{n}{k} \rfloor \} 
    \Big\} .
\end{eqnarray*}

\noindent Recall that $h_m - h_{m-1} = h_1$. Applying the strong Markov property successively at times $\mathscrHS_{h_{k-1}}$, $\mathscrHS_{h_{k-2}}$, $\cdots$, $\mathscrHS_{h_1}$, this gives that\footnote{For the term $m=k$ on the right-hand side, there is no need to consider $\{ \Delta S_{\mathscrHS_x} \le r^\theta\}$, whereas the $m=1$ term has only the value $x=h_1$. The current form of the inequality is used to give a compact expression for the lower bound.}
\begin{eqnarray}
    q_1(r)
 &\ge&\prod_{m=1}^k \inf_{x\in (r^\theta, \, h_1]} 
    \Q \Big\{ 
    \bigcap_{j=0}^{\mathscrHS_x-1} 
    \{ \overline{S}_j - S_j < \lambda_m, \, S_j \ge - \beta\} \cap
    \nonumber
    \\
 &&\qquad\qquad  
    \cap \{ \Delta S_{\mathscrHS_x} \le r^\theta, \, 
    \mathscrHS_x < \lfloor \frac{n}{k} \rfloor \} \Big\} \, .
    \label{q1(r)>}
\end{eqnarray}


\noindent We let $r\to \infty$. By Lemma \ref{l:RW-chute-lb}, uniformly in $m\in [1, \, k]$ and $x\in (r^\theta, \, h_1]$,
\begin{eqnarray*}
    \Q \Big\{ 
    \bigcap_{j=0}^{\mathscrHS_x-1} 
    \{ \overline{S}_j - S_j < \lambda_m, \, S_j \ge - \beta\}
    \Big\}
 &\ge&\exp\Big[ -(1+o(1)) \frac{x}{\lambda_m} \Big]
    \\
 &\ge& \exp\Big[ -(1+o(1)) \frac{r^\chi}{\lambda_m} \Big] .
\end{eqnarray*}


\noindent On the other hand, (\ref{moment-expo-overshoot}) tells us that $c_{17} := \sup_{b>0} \E_\Q[ \exp( c_2 \, \Delta S_{\mathscrHS_b})] <\infty$. By the Markov inequality, for $r\to \infty$, uniformly in $m\in [1, \, k]$ and $x\in (r^\theta, \, h_1]$,
$$
    \Q \{ \Delta S_{\mathscrHS_x} > r^\theta \}
    \le
    c_{17} \, \ee^{- c_2 \, r^\theta}
    \le 
    \frac13 \exp\Big[ -(1+o(1)) \frac{r^\chi}{\lambda_m} \Big] .
$$

\noindent [The last inequality, valid for all sufficiently large $r$, relies on the facts that $\theta > \frac{\chi}{2}$ and that $\lambda_m \ge (2r^\chi)^{1/2}$.] Also, for some constant $c_{18}>0$ and all sufficiently large $r$ and all $m\in [1, \, k]$, $\sup_{x\in (r^\theta, \, h_1]}\Q\{ \mathscrHS_x \ge \lfloor \frac{n}{k} \rfloor \} \le c_{18} \, \frac{h_1}{(\lfloor \frac{n}{k} \rfloor)^{1/2}}$ (see Theorem A of Kozlov~\cite{kozlov}), which is bounded by $\frac13 \exp [ -(1+o(1)) \frac{r^\chi}{\lambda_m} ]$ as well for some constant $\varepsilon_1>0$ (for $r\to \infty$; recalling that $n:= \lfloor \ee^{\varepsilon_1\, r^{1/2}} \rfloor$). [We use the fact that $\frac12 > \frac{\chi}{2}$.] As a consequence, for $r\to \infty$, uniformly in $m\in [1, \, k]$ and $x\in (r^\theta, \, h_1]$,
\begin{eqnarray*}
 &&\Q \Big\{ 
    \bigcap_{j=0}^{\mathscrHS_x-1} 
    \{ \overline{S}_j - S_j < \lambda_m, \, S_j \ge - \beta\} \cap
    \{ \Delta S_{\mathscrHS_x} \le r^\theta, \, 
    \mathscrHS_x < \lfloor \frac{n}{k} \rfloor \} \Big\} 
    \\
 &\ge& \Q \Big\{ 
    \bigcap_{j=0}^{\mathscrHS_x-1} 
    \{ \overline{S}_j - S_j < \lambda_m, \, S_j \ge - \beta\} \Big\} 
    - 
    \Q\{ \Delta S_{\mathscrHS_x} > r^\theta \}
    - 
    \Q\{ \mathscrHS_x \ge \lfloor \frac{n}{k} \rfloor \} 
    \\
 &\ge&\frac13 \exp\Big[ -(1+o(1)) \frac{r^\chi}{\lambda_m} \Big] , 
\end{eqnarray*}

\noindent which is still $\exp[ -(1+o(1)) \frac{r^\chi}{\lambda_m}]$ by changing the value of $o(1)$. Going back to (\ref{q1(r)>}), we see that for $r\to \infty$,
\begin{equation}
    q_1(r)
    \ge
    \exp\Big[ -(1+o(1)) \sum_{m=1}^k \frac{r^\chi}{\lambda_m} \Big]
    =
    \ee^{-(1+o(1))(2r)^{1/2}},
    \label{q1(r)>:2}
\end{equation}

\noindent the last identity following from the observation in (\ref{sum(hm-hm-1)/lambda_m=sqrt(2r)}) that $\sum_{m=1}^k \frac{r^\chi}{\lambda_m} \sim (2r)^{1/2}$, $r\to \infty$. 

We now estimate $q_2(r)$. By definition,
\begin{eqnarray}
    q_2(r)
 &\le&\sum_{\ell=1}^n \Q \Big[ 
    \bigcap_{j=0}^{\mathscrHS_r-1} \{ \overline{S}_j - S_j < a^{(S)}_j\}; 
    \;
    \max_{1\le i<k} \Delta S_{\mathscrHS_{h_i}} \le r^\theta ;
    \; 
    \eta_\ell > \ee^{\varepsilon r^{1/2}} ;
    \;
    \ell \le \mathscrHS_r \Big]
    \nonumber
    \\
 &=&\sum_{\ell=1}^n \sum_{m=1}^k q_2^{(\ell, \, m)}(r) \, ,
    \label{q2<q3}
\end{eqnarray}

\noindent where
\begin{eqnarray*}
 &&q_2^{(\ell, \, m)}(r)
    \\
 &:=& \Q \Big[ 
    \bigcap_{j=0}^{\mathscrHS_r-1} \{ \overline{S}_j - S_j < a^{(S)}_j\} ;
    \;  
    \max_{1\le i<k} \Delta S_{\mathscrHS_{h_i}} \le r^\theta; 
    \;
    \eta_\ell > \ee^{\varepsilon r^{1/2}} ;
    \;
    \mathscrHS_{h_{m-1}} < \ell \le \mathscrHS_{h_m}
    \Big]
    \\
 &=&\Q \Big[ 
    \bigcap_{i=1}^k \bigcap_{j=\mathscrHS_{h_{i-1}}}^{\mathscrHS_{h_i}-1} \{ \overline{S}_j - S_j < \lambda_i\} ;
    \;  
    \max_{1\le i<k} \Delta S_{\mathscrHS_{h_i}} \le r^\theta; 
    \;
    \eta_\ell > \ee^{\varepsilon r^{1/2}} ;
    \;
    \mathscrHS_{h_{m-1}} < \ell \le \mathscrHS_{h_m}
    \Big] .
\end{eqnarray*}

\noindent We apply the strong Markov property at $\mathscrHS_{h_{k-1}}$, to see that for $1\le m<k$,
\begin{eqnarray*}
    q_2^{(\ell, \, m)}(r)
 &\le& \Q\Big[ 
    \bigcap_{i=1}^{k-1} \bigcap_{j=\mathscrHS_{h_{i-1}}}^{\mathscrHS_{h_i}-1} \{ \overline{S}_j - S_j < \lambda_i\} ;
    \;  
    \max_{1\le i<k} \Delta S_{\mathscrHS_{h_i}} \le r^\theta; 
    \;
    \eta_\ell > \ee^{\varepsilon r^{1/2}} ;
    \;
    \\
 &&\qquad
    \mathscrHS_{h_{m-1}} < \ell \le \mathscrHS_{h_m}
    \Big]
    \times
    \sup_{x\in [h_{k-1}, \, h_{k-1} + r^\theta]} 
    \Q \Big[ \bigcap_{j=0}^{\mathscrHS_{h_k-x}-1} \{ \overline{S}_j - S_j < \lambda_k\} \Big].
\end{eqnarray*}

\noindent Let $r\to \infty$. By Lemma \ref{l:RW-chute-ub}, we have, uniformly in $x\in [h_{k-1}, \, h_{k-1} + r^\theta]$,
\begin{eqnarray*}
    \Q \Big[ \bigcap_{j=0}^{\mathscrHS_{h_k-x}-1} \{ \overline{S}_j - S_j < \lambda_k\} \Big]
 &\le&\exp \Big[ - (1+o(1))\, \frac{h_k-h_{k-1} - r^\theta}{\lambda_k} \Big]
    \\
 &\le& \exp \Big[ - (1+o(1))\, \frac{r^\chi}{\lambda_k} \Big] \, .
\end{eqnarray*}

\noindent We iterate the argument and apply the strong Markov property successively at $\mathscrHS_{h_{k-2}}$, $\mathscrHS_{h_{k-3}}$, $\cdots$, $\mathscrHS_{h_m}$, to see that
\begin{eqnarray*}
    q_2^{(\ell, \, m)}(r)
 &\le& \Q\Big[ 
    \bigcap_{i=1}^m \bigcap_{j=\mathscrHS_{h_{i-1}}}^{\mathscrHS_{h_i}-1} \{ \overline{S}_j - S_j < \lambda_i\} ;
    \;  
    \max_{1\le i\le m} \Delta S_{\mathscrHS_{h_i}} \le r^\theta; 
    \;
    \eta_\ell > \ee^{\varepsilon r^{1/2}} ;
    \;
    \\
 &&\qquad
    \mathscrHS_{h_{m-1}} < \ell \le \mathscrHS_{h_m}
    \Big]
    \times
    \exp \Big[ - (1+o(1)) \sum_{i=m+1}^k \frac{r^\chi}{\lambda_i} \Big]
    \\
 &\le& \Q\Big[ 
    \bigcap_{i=1}^{m-1} \bigcap_{j=\mathscrHS_{h_{i-1}}}^{\mathscrHS_{h_i}-1} \{ \overline{S}_j - S_j < \lambda_i\} ;
    \;  
    \max_{1\le i\le m-2} \Delta S_{\mathscrHS_{h_i}} \le r^\theta; 
    \;
    \eta_\ell > \ee^{\varepsilon r^{1/2}} ;
    \;
    \\
 &&\qquad
    \mathscrHS_{h_{m-1}} < \ell \Big]
    \times
    \exp \Big[ - (1+o(1)) \sum_{i=m+1}^k \frac{r^\chi}{\lambda_i} \Big] .
\end{eqnarray*}

\noindent To bound the probability expression $\Q[ \cdots]$ on the right-hand side, we note that under $\Q$, given $\mathscrHS_{h_{m-1}} < \ell$, $\eta_\ell$ is independent of everything concerning the potential $V(\cdot)$ until $\mathscrHS_{h_{m-1}}$, and has the law of $\eta_1$. Consequently,
\begin{eqnarray*}
    q_2^{(\ell, \, m)}(r)
 &\le& \Q\Big[ 
    \bigcap_{i=1}^{m-1} \bigcap_{j=\mathscrHS_{h_{i-1}}}^{\mathscrHS_{h_i}-1} \{ \overline{S}_j - S_j < \lambda_i\} ;
    \;  
    \max_{1\le i\le m-2} \Delta S_{\mathscrHS_{h_i}} \le r^\theta; 
    \;
    \mathscrHS_{h_{m-1}} < \ell \Big]
    \times
    \\
 &&\qquad
    \times \Q(\eta_1 > \ee^{\varepsilon r^{1/2}})
    \times
    \exp \Big[ - (1+o(1)) \sum_{i=m+1}^k \frac{r^\chi}{\lambda_i} \Big] 
    \\
 &\le&\Q\Big[ 
    \bigcap_{i=1}^{m-1} \bigcap_{j=\mathscrHS_{h_{i-1}}}^{\mathscrHS_{h_i}-1} \{ \overline{S}_j - S_j < \lambda_i\} ;
    \;  
    \max_{1\le i\le m-2} \Delta S_{\mathscrHS_{h_i}} \le r^\theta \Big]
    \times
    \\ 
 &&\qquad
    \times \Q(\eta_1 > \ee^{\varepsilon r^{1/2}})
    \times
    \exp \Big[ - (1+o(1)) \sum_{i=m+1}^k \frac{r^\chi}{\lambda_i} \Big] .
\end{eqnarray*}


Looking at the two probability expressions $\Q[ \cap_{i=1}^{m-1} \cdots ]$ and $\Q(\eta_1 > \ee^{\varepsilon r^{1/2}})$ on the right-hand side. The first probability expression is,
according to (\ref{f<:m=0}), 
bounded by $\exp [ - (1+o(1)) \sum_{\ell=1}^{m-1} \frac{r^\chi}{\lambda_\ell} ]$. For the second probability expression, let us recall that $\eta_1 = \sum_{y: \, |y|=1} \ee^{-V(y)}$ by definition; so by (\ref{nu}), there exists a constant $c_{19}>0$ such that $\Q(\eta_1 > \ee^{\varepsilon r^{1/2}}) \le c_{19} \, \ee^{-c_1 \, \varepsilon r^{1/2}}$. We have thus proved that, for $1\le m\le k$,
\begin{eqnarray*}
    q_2^{(\ell, \, m)}(r)
 &\le& c_{19} \, \ee^{-c_1 \, \varepsilon r^{1/2}} \exp \Big[ - (1+o(1)) \sum_{i: \, 1\le i\le k, \, i\not= m} \frac{r^\chi}{\lambda_i} \Big] 
    \\
 &\le& c_{19} \, \ee^{-c_1 \, \varepsilon r^{1/2} - (1+o(1)) (2r)^{1/2} } . 
\end{eqnarray*}

\noindent Since $q_2(r) \le \sum_{\ell=1}^n \sum_{m=1}^k q_2^{(\ell, \, m)}(r)$ (see (\ref{q2<q3})), and $n:= \lfloor \ee^{\varepsilon_1 \, r^{1/2}}\rfloor \le \ee^{\varepsilon_1 \, r^{1/2}}$, this yields that
$$
q_2(r)
\le
c_{19} \, k\, \ee^{-(c_1 \, \varepsilon - \varepsilon_1) r^{1/2} - (1+o(1)) (2r)^{1/2} } \, .
$$

\noindent Recall that $\E[ E_\omega(Z_r) ] \ge \frac{q_1(r) - q_2(r)}{n}$ (see (\ref{E(Z)>q1-q2})) and that $q_1(r) \ge \ee^{-(1+o(1))(2r)^{1/2}}$ (see (\ref{q1(r)>:2})), we obtain that for $r\to \infty$,
$$
\E[ E_\omega(Z_r) 
    ]
\ge
\frac{1}{n} \Big[ \ee^{-(1+o(1))(2r)^{1/2}} - c_{19} \, k \, \ee^{-(c_1 \, \varepsilon - \varepsilon_1) r^{1/2} - (1+o(1)) (2r)^{1/2} } \Big].
$$

\noindent Since $\varepsilon_1 \in (0, \, c_1 \, \varepsilon)$, the term $c_{19} \, k \, \ee^{-(c_1 \, \varepsilon - \varepsilon_1) r^{1/2} - (1+o(1)) (2r)^{1/2} }$ does not play any role when taking the limit $r\to \infty$ (recalling that $k:= \lfloor r^{1-\chi}\rfloor$). By definition, $n:= \lfloor \ee^{\varepsilon_1 \, r^{1/2}}\rfloor$, this readily yields (\ref{E(Z)}).\hfill$\Box$

\subsection{Proof of Lemma \ref{l:Z}: the inequality (\ref{E(Z2)})}
\label{subs:proof-second-lemma-(ii)}

Recall definition again from (\ref{Z}): $Z_r := \sum_{x\in \mathscrHV_r^*} {\bf 1}_{\{ T_x < T_{{\buildrel \leftarrow \over \varnothing}} \} }$, where 
\begin{eqnarray*}
    \mathscrHV_r^*
 &:=& \Big\{ x\in \mathscrHV_r: \, \max_{1\le m<k} \Delta V(x_{\mathscrHx_{h_m}}) \le r^\theta , \; \underline{V}(x) \ge - \beta , \; |x| < \lfloor \ee^{\varepsilon_1 \, r^{1/2}}\rfloor \, , \;
    \\
 && \overline{V}(x_j) - V(x_j) \le a_j^{(x)} , \; \forall 0\le j<|x| , \;  \max_{0\le j<|x|} \Lambda(x_j) \le \ee^{\varepsilon r^{1/2}} \Big\} ,
\end{eqnarray*}

\noindent with $\Lambda(x) := \sum_{y: \, {\buildrel \leftarrow \over y} = x} \ee^{-\Delta V(y)}$ as in (\ref{Lambda}). By definition,
\begin{eqnarray}
    E_\omega  (Z_r^2)
 &=&\sum_{x, \, y \in \mathscrHV_r^*} 
    P_\omega\{ T_x < T_{{\buildrel \leftarrow \over \varnothing}}, \, T_y < T_{{\buildrel \leftarrow \over \varnothing}} \} 
    \nonumber
    \\
 &=& E_\omega(Z_r) + \sum_{x\not= y \in \mathscrHV_r^*} 
     P_\omega\{ T_x < T_{{\buildrel \leftarrow \over \varnothing}}, \, T_y < T_{{\buildrel \leftarrow \over \varnothing}} \} \, .
     \label{E(Z2)<:start}
\end{eqnarray}

By (\ref{P(T_x<T_0)=}), $P_\omega\{ T_x < T_{{\buildrel \leftarrow \over \varnothing}}\} \le \ee^{-\overline{V}(x)}$. On the other hand, by the definition of $\mathscrHV_r$, we have $\overline{V}(x) = V(x)$ for $x\in \mathscrHV_r^* \subset \mathscrHV_r$. So
\begin{eqnarray*}
    E_\omega(Z_r)
 &\le&\sum_{x\in \mathscrHV_r^*} \ee^{- V(x)}
    \\
 &\le&\sum_{x\in \mathscrHV_r} \ee^{- V(x)} \, 
    {\bf 1}_{\{ \max_{1\le m<k} \Delta V(x_{\mathscrHx_{h_m}}) \le r^\theta\} } \,
    {\bf 1}_{\{ \overline{V}(x_j) - V(x_j) \le a_j^{(x)} , \; \forall 0\le j<|x| \} } \, .
\end{eqnarray*}

\noindent Taking expectation on both sides, we obtain that 
$$
\E[E_\omega(Z_r)]
\le
\E \Big( \sum_{x\in \mathscrHV_r} \ee^{- V(x)} \, {\bf 1}_{\{ \max_{1\le m<k} \Delta V(x_{\mathscrHx_{h_m}}) \le r^\theta\} } \,
    {\bf 1}_{\{ \overline{V}(x_j) - V(x_j) \le a_j^{(x)} , \; \forall 0\le j<|x| \} } \Big),
$$

\noindent which, by formula (\ref{many-to-one:stopping-lines}), is
$$
=
\Q \Big( \max_{1\le m<k} \Delta S_{\mathscrHS_{h_m}} \, , \overline{S}_j - S_j \le a_j^{(S)} , \; \forall 0\le j< \mathscrHS_r \Big) .
$$

\noindent Applying (\ref{f<:m=0}), we get that $\E[E_\omega(Z_r)] \le \ee^{- (1+o(1)) \sum_{\ell=1}^k \frac{r^\chi}{\lambda_\ell}}$. Since $\sum_{\ell=1}^k \frac{r^\chi}{\lambda_\ell} \sim (2r)^{1/2}$ (see (\ref{sum(hm-hm-1)/lambda_m=sqrt(2r)})), we arrive at:
\begin{equation}
    \E[E_\omega(Z_r)]
    \le
    \ee^{-(1+o(1))(2r)^{1/2}} \, .
    \label{E(Z)<}
\end{equation}

Also, since $V(x) \ge r$ for $x\in \mathscrHV_r^*$, we have $\sum_{x\in \mathscrHV_r^*} \ee^{-2V(x)} \le \ee^{-r} \sum_{x\in \mathscrHV_r^*} \ee^{- V(x)}$, so for all sufficiently large $r$,
\begin{equation}
    \E\Big( \sum_{x\in \mathscrHV_r^*} \ee^{-2V(x)} \Big)
    \le
    \ee^{-r} \, .
    \label{E(second-moment-diagonal)}
\end{equation}

\noindent By (\ref{E(Z)<}) and (\ref{E(Z2)<:start}), we have
\begin{equation}
    \E[E_\omega  (Z_r^2)]
    \le
    \ee^{-(1+o(1))(2r)^{1/2}} 
    +
    \E \Big[ 
    \sum_{x\not= y \in \mathscrHV_r^*} 
     P_\omega\{ T_x < T_{{\buildrel \leftarrow \over \varnothing}}, \, T_y < T_{{\buildrel \leftarrow \over \varnothing}} \} \Big] \, .
    \label{E(Z2)<}
\end{equation}

For any pair of distinct vertices $x\not= y$, let $x\wedge y$ denote their youngest common ancestor; equivalently, $x\wedge y$ is the unique vertex satisfying  $[\![ \varnothing , \, x\wedge y ]\!] = [\![ \varnothing , \, x ]\!] \cap [\![ \varnothing , \, y ]\!]$. Consider the quenched probability expression
$$
P_\omega\{ T_x < T_y < T_{{\buildrel \leftarrow \over \varnothing}} \} .
$$

\noindent To realize $T_x < T_y < T_{{\buildrel \leftarrow \over \varnothing}}$, the biased walk first needs to hit $x\wedge y$ before hitting ${\buildrel \leftarrow \over \varnothing}$, then, starting from $x\wedge y$, it should hit $x$ before hitting ${\buildrel \leftarrow \over \varnothing}$, (and then, starting from $x$, it hits automatically $x\wedge y$ before hitting ${\buildrel \leftarrow \over \varnothing}$), and then, starting from $x\wedge y$, it should hit $y$ before hitting ${\buildrel \leftarrow \over \varnothing}$. Applying the strong Markov property, we obtain that
$$
P_\omega\{ T_x < T_y < T_{{\buildrel \leftarrow \over \varnothing}} \}
\le
P_\omega\{ T_{x\wedge y} < T_{{\buildrel \leftarrow \over \varnothing}} \}
P_\omega^{x\wedge y}\{ T_x < T_{{\buildrel \leftarrow \over \varnothing}} \}
P_\omega^{x\wedge y}\{ T_y < T_{{\buildrel \leftarrow \over \varnothing}} \} ,
$$

\noindent where, for any vertex $z$, $P_\omega^z$ denotes the (quenched) probability under which the biased walk starts at $z$. By exchanging $x$ and $y$, we also have
$$
P_\omega\{ T_y < T_x < T_{{\buildrel \leftarrow \over \varnothing}} \}
\le
P_\omega\{ T_{x\wedge y} < T_{{\buildrel \leftarrow \over \varnothing}} \}
P_\omega^{x\wedge y}\{ T_y < T_{{\buildrel \leftarrow \over \varnothing}} \}
P_\omega^{x\wedge y}\{ T_x < T_{{\buildrel \leftarrow \over \varnothing}} \} .
$$

\noindent Hence
\begin{eqnarray*}
    P_\omega\{ T_x < T_{{\buildrel \leftarrow \over \varnothing}}, \, T_y < T_{{\buildrel \leftarrow \over \varnothing}} \}
 &=& P_\omega\{ T_x < T_y < T_{{\buildrel \leftarrow \over \varnothing}} \} +
    P_\omega\{ T_y < T_x < T_{{\buildrel \leftarrow \over \varnothing}} \}
    \\
 &\le&2P_\omega\{ T_{x\wedge y} < T_{{\buildrel \leftarrow \over \varnothing}} \}
P_\omega^{x\wedge y}\{ T_x < T_{{\buildrel \leftarrow \over \varnothing}} \}
P_\omega^{x\wedge y}\{ T_y < T_{{\buildrel \leftarrow \over \varnothing}} \} .
\end{eqnarray*}

\noindent [Although we have implicitly assumed ${x\wedge y}$ is different from the root $\varnothing$, the last inequality remains trivially valid even if ${x\wedge y}$ is the root.] By (\ref{P(T_x<T_0)=}), $P_\omega\{ T_{x\wedge y} < T_{{\buildrel \leftarrow \over \varnothing}} \} \le \ee^{-\overline{V}(x\wedge y)}$. More generally, we use (\ref{P(T_x<T_0):cas_general}) to see that 
$$
P_\omega^{x\wedge y}\{ T_x < T_{{\buildrel \leftarrow \over \varnothing}} \}
\le
(|x\wedge y|+1)\ee^{-[\overline{V}(x)- \overline{V}(x\wedge y)]} \, .
$$

\noindent We also have $P_\omega^{x\wedge y}\{ T_y < T_{{\buildrel \leftarrow \over \varnothing}} \} \le (|x\wedge y|+1) \ee^{-[\overline{V}(y)- \overline{V}(x\wedge y)]}$ by interchanging the roles of $x$ and $y$. As a consequence,
$$
P_\omega\{ T_x < T_{{\buildrel \leftarrow \over \varnothing}}, \, T_y < T_{{\buildrel \leftarrow \over \varnothing}} \}
\le
2(|x\wedge y|+1)^2\, \ee^{\overline{V}(x\wedge y)-\overline{V}(x)-\overline{V}(y)} ,
$$

\noindent which is bounded by $2(|x\wedge y|+1)^2 \,\ee^{\overline{V}(x\wedge y)- V(x) -V(y)}$. Moreover, for $x\in \mathscrHV_r^*$, we have $|x\wedge y| +1 \le |x| +1 \le \lfloor \ee^{\varepsilon_1 \, r^{1/2}}\rfloor$. Going back to (\ref{E(Z2)<}), we obtain that
\begin{eqnarray}
 &&\E[E_\omega  (Z_r^2)]
    \nonumber
    \\
 &\le&\ee^{-(1+o(1))(2r)^{1/2}}
    +
    2\ee^{2\varepsilon_1 \, r^{1/2}}\, 
    \E \Big( \sum_{z: \; \overline{V}(z) < r} 
    \sum_{x, \, y \in \mathscrHV_r^*: \; x\wedge y = z} 
    \ee^{\overline{V}(z)-V(x)-V(y)} \Big)
    \label{Zr:second_moment:preliminaire}
    \\
 &=& \ee^{-(1+o(1))(2r)^{1/2}}
    +
    2\ee^{2\varepsilon_1 \, r^{1/2}}\, 
    \E\Big( \sum_{n=0}^\infty \sum_{m=1}^k 
    \Sigma_3^{(n,m)} \Big) ,
    \label{Zr:second_moment}
\end{eqnarray}

\noindent where
$$
\Sigma_3^{(n,m)}
:=
\sum_{z: \; |z|=n}
\ee^{\overline{V}(z)} \,
{\bf 1}_{\{ h_{m-1} \le \overline{V}(z) < h_m\} } 
\sum_{x, \, y \in \mathscrHV_r^*: \; x\wedge y = z} 
\ee^{-V(x)-V(y)} \, .
$$

\noindent For further use, we also see from the inequality $E_\omega(Z_r) \le \sum_{x\in \mathscrHV_r^*} \ee^{- V(x)}$ that, for all sufficiently large $r$,
\begin{equation}
    \E[(E_\omega Z_r)^2]
    \le
    \ee^{-r}
    +
    \E\Big( \sum_{z: \; \overline{V}(z) < r} 
{\bf 1}_{\{ \underline{V}(z) \ge -\beta\} }
\sum_{x, \, y \in \mathscrHV_r^*: \; x\wedge y = z} 
\ee^{-V(x)-V(y)} \Big) .
    \label{Zr:second_moment-bis}
\end{equation}

\noindent The term $\ee^{-r}$ comes from $\E ( \sum_{x\in \mathscrHV_r^*} \ee^{-2V(x)})$ and (\ref{E(second-moment-diagonal)}). The indicator function ${\bf 1}_{\{ \underline{V}(z) \ge -\beta\} }$ was implicitly present in $x \in \mathscrHV_r^*$; it is written explicitly here because it is going to play a crucial role later. We note that the expectation expressions on the right-hand side of (\ref{Zr:second_moment:preliminaire}) and (\ref{Zr:second_moment-bis}) are very similar to each other, except that there is no $\overline{V}(z)$ term on the right-hand side of (\ref{Zr:second_moment-bis}).

For each pair $(n, \, m)$, we estimate $\E(\Sigma_3^{(n,m)})$. By definition (recalling that $x_i$ is the ancestor of $x$ in generation $i$ for $i\le |x|$),
\begin{eqnarray*}
    \Sigma_3^{(n,m)}
 &=&\sum_{z: \; |z|=n} 
    \ee^{\overline{V}(z)}\,
    {\bf 1}_{\{ h_{m-1} \le \overline{V}(z) < h_m\} } 
    \sum_{u\not= v, \, {\buildrel \leftarrow \over u}=z={\buildrel \leftarrow \over v} } \ee^{-V(u)-V(v)} \times
    \\
 &&\qquad \times
    \sum_{x\in \mathscrHV_r^*: \; x_{n+1}=u} \ee^{-[V(x)-V(u)]} 
    \sum_{y\in \mathscrHV_r^*: \; y_{n+1}=v} \ee^{-[V(y)-V(v)]} \, .
\end{eqnarray*}

\noindent We first take expectation conditioning on $\mathscr{F}_{n+1} := \sigma \{ V(w): \, |w| \le n+1\}$, the $\sigma$-field generated by the random potential in the first $n+1$ generations:
\begin{eqnarray}
 &&\E (\Sigma_3^{(n,m)} \, | \, \mathscr{F}_{n+1})
    \nonumber
    \\
 &\le& \sum_{z: \; |z|=n} 
    \ee^{\overline{V}(z)} \,
    {\bf 1}_{\{ h_{m-1} \le \overline{V}(z) < h_m\} } \,
    {\bf 1}_{\{ \overline{V}(z_i) - V(z_i) < a_i^{(z)}, \; \forall 0\le i\le n\} } \,
    {\bf 1}_{\{ \max_{1\le \ell < m} \Delta V(z_{\mathscrHz_{h_\ell}}) \le r^\theta\} }
    \times
    \nonumber
    \\
 &&\qquad \times
    {\bf 1}_{\{ \Lambda(z) \le \ee^{\varepsilon r^{1/2}} \} }\,
    \sum_{(u, \, v): \, u\not= v, \, {\buildrel \leftarrow \over u}=z={\buildrel \leftarrow \over v} } \ee^{-V(u)-V(v)} \,
    f_m(V(u)) f_m(V(v)),
    \label{Sigma_3}
\end{eqnarray}

\noindent where $\Lambda(x) := \sum_{y: \, {\buildrel \leftarrow \over y} = x} \ee^{-\Delta V(y)}$ as in (\ref{Lambda}), and for $s<h_{m+1}$,
$$
f_m(s)
:=
\E \Big\{ 
\sum_{x\in \mathscrHV_{r-s}} \ee^{-V(x)} 
\Big( \prod_{\ell=m+1}^{k-1}{\bf 1}_{\{ \Delta V(x_{\mathscrHx_{h_\ell-s}}) \le r^\theta\} } \Big)
\Big( \prod_{\ell=m+2}^k \prod_{i=\mathscrHx_{h_{\ell-1}-s}}^{\mathscrHx_{h_\ell-s} -1} {\bf 1}_{ \{ \overline{V}(x_i) - V(x_i) <\lambda_\ell \} } \Big) \Big\} . 
$$

\noindent Some care needs to be taken in order to make (\ref{Sigma_3}) valid in all situations. On the right-hand side of (\ref{Sigma_3}), $V(u) < h_m$ for {\bf most} $u$ with ${\buildrel \leftarrow \over u}=z$ (and $V(u)<r$ for most $v$ with ${\buildrel \leftarrow \over v}=z$); however, there is a possible situation when $V(u) \ge h_m$: this is when $u\in \mathscrHV_{h_m}$ (for some $1\le m\le k$), in which case we only have $V(u) \le h_m + r^\theta$ (which is strictly smaller than $h_{m+1}$). In order to take care of this situation, only overshoots $\Delta V(x_{\mathscrHx_{h_\ell-s}})$ for $\ell >m$ are involved in the definition of $f_m(s)$. In particular, $f_{k-1}(s) = 1$ for $s<r$, and $f_k(s)$ should be defined as $1$ for all $s\in \r$.

By the formula (\ref{many-to-one:stopping-lines}), this gives that for $s<h_{m+1}$,
$$
f_m(s)
=
\Q\Big( \bigcap_{\ell=m+1}^{k-1}{\bf 1}_{\{ \Delta S_{\mathscrHS_{h_\ell-s}} \le r^\theta\} } 
\cap
\bigcap_{\ell=m+2}^k \bigcap_{i=\mathscrHS_{h_{\ell-1}-s}}^{\mathscrHS_{h_\ell-s} -1} \{ \overline{S}_i -S_i <\lambda_\ell \}\Big) ,
$$

\noindent where $\mathscrHS_t :=  \inf\{ i\ge 0: \, S_i \ge t\}$ (for any $t\ge 0$) as in (\ref{mathscrHS_r}). By Claim \ref{claim:proba}, we arrive at the following estimate: when $r\to \infty$,
$$
f_m(s)
\le
\exp\Big( - (1+o(1)) \sum_{\ell=m+2}^k \frac{r^\chi}{\lambda_\ell} \Big) ,
$$

\noindent uniformly in $s<h_{m+1}$ and $m\in [1, \, k]$ (and in $n\ge 1$).

Let us go back to (\ref{Sigma_3}), and first look at the double sum $\sum_{(u, \, v): \, u\not= v, \, {\buildrel \leftarrow \over u}=z={\buildrel \leftarrow \over v} }$ on the right-hand side. Thanks to the upper bound for $f_m(s)$ we have just obtained that is valid uniformly in $s\ge 0$, we get that, on the right-hand side of (\ref{Sigma_3}),
\begin{eqnarray*}
 &&{\bf 1}_{\{ \Lambda(z) \le \ee^{\varepsilon r^{1/2}} \} }\, 
    \sum_{(u, \, v): \, u\not= v, \, {\buildrel \leftarrow \over u}=z={\buildrel \leftarrow \over v} } \ee^{-V(u)-V(v)} \,
    f_m(V(u)) f_m(V(v))
    \\
 &\le&{\bf 1}_{\{ \Lambda(z) \le \ee^{\varepsilon r^{1/2}} \} }\,
    \ee^{- (2+o(1)) \sum_{\ell=m+2}^k \frac{r^\chi}{\lambda_\ell} }
    \Big[ \sum_{u: \, {\buildrel \leftarrow \over u}=z} \ee^{-V(u)} \Big]^2
    \\
 &\le& \ee^{- (2+o(1)) \sum_{\ell=m+2}^k \frac{r^\chi}{\lambda_\ell} } 
    \Big[ \ee^{-V(z)} \, \ee^{\varepsilon r^{1/2}} \Big]^2 \, ,
\end{eqnarray*}

\noindent where, in the last inequality, we used the definition of $\Lambda(z):= \sum_{u: \, {\buildrel \leftarrow \over u}=z} \ee^{-[V(u)-V(z)]}$ as in (\ref{Lambda}) to see that on the event $\{ \Lambda(z) \le \ee^{\varepsilon r^{1/2}} \}$, we have $\sum_{u: \, {\buildrel \leftarrow \over u}=z} \ee^{-V(u)} = \ee^{-V(z)} \, \Lambda(z) \le \ee^{-V(z)} \, \ee^{\varepsilon r^{1/2}}$. Therefore, (\ref{Sigma_3}) yields that
\begin{eqnarray*}
    \E (\Sigma_3^{(n,m)} \, | \, \mathscr{F}_{n+1})
 &\le& \ee^{2\varepsilon r^{1/2}- (2+o(1)) \sum_{\ell=m+2}^k \frac{r^\chi}{\lambda_\ell} }
    \sum_{z: \; |z|=n} 
    \ee^{\overline{V}(z)-2V(z)} \,
    {\bf 1}_{\{ h_{m-1} \le \overline{V}(z) < h_m\} } \,
    \times
    \\
 && \times 
    {\bf 1}_{\{ \overline{V}(z_i) - V(z_i) < a_i^{(z)}, \; \forall 0\le i\le n\} } \,
    {\bf 1}_{\{ \max_{1\le \ell < m} \Delta V(z_{\mathscrHz_{h_\ell}}) \le r^\theta\} } \, .
\end{eqnarray*}

\noindent Taking expectation to get rid of the conditioning, and using the many-to-one formula (\ref{many-to-one}), we obtain that
\begin{eqnarray*}
    \E (\Sigma_3^{(n,m)})
 &\le& \ee^{2\varepsilon r^{1/2}- (2+o(1)) \sum_{\ell=m+2}^k \frac{r^\chi}{\lambda_\ell} } \,
    \E_{\Q} \Big[ 
    \ee^{\overline{S}_n -S_n} \,
    {\bf 1}_{\{ h_{m-1} \le \overline{S}_n < h_m\} } \,
    \times
    \\
 && \times 
    {\bf 1}_{\{ \overline{S}_i - S_i < a_i^{(S)}, \; \forall 0\le i\le n\} } \,
    {\bf 1}_{\{ \max_{1\le \ell < m} \Delta S_{\mathscrHS_{h_\ell}} \le r^\theta\} } 
    \Big] .
\end{eqnarray*}

\noindent Going back to (\ref{Zr:second_moment}), this yields that
\begin{eqnarray}
 &&\E[E_\omega  (Z_r^2)]
    \le
    \ee^{-(1+o(1))(2r)^{1/2}}
    +
    2\ee^{2\varepsilon_1 \, r^{1/2}}\, 
    \sum_{n=0}^\infty \sum_{m=1}^k 
    \ee^{2\varepsilon r^{1/2}- (2+o(1)) \sum_{\ell=m+2}^k \frac{r^\chi}{\lambda_\ell} } \times
    \nonumber
    \\
 &&\qquad \times 
    \E_{\Q} \Big[ 
    \ee^{\overline{S}_n -S_n} \,
    {\bf 1}_{\{ h_{m-1} \le \overline{S}_n < h_m\} } \, 
    {\bf 1}_{\{ \overline{S}_i - S_i < a_i^{(S)}, \; \forall 0\le i\le n\} } \,
    {\bf 1}_{\{ \max_{1\le \ell < m} \Delta S_{\mathscrHS_{h_\ell}} \le r^\theta\} } 
    \Big] .
    \label{Zr:second_moment-2}
\end{eqnarray}

\noindent Similarly, (\ref{Zr:second_moment-bis}) leads to: for $r\to \infty$,
\begin{eqnarray}
 &&\E[(E_\omega Z_r)^2]
    \le
    \ee^{-r}
    +
    \sum_{n=0}^\infty \sum_{m=1}^k 
    \ee^{2\varepsilon r^{1/2}- (2+o(1)) \sum_{\ell=m+2}^k \frac{r^\chi}{\lambda_\ell} } \,
    \E_{\Q} \Big[ 
    \ee^{-S_n} \,
    {\bf 1}_{\{ \min_{0\le i\le n} S_i \ge -\beta\} } \, \times
    \nonumber
    \\
 &&\qquad\qquad \times 
    {\bf 1}_{\{ h_{m-1} \le \overline{S}_n < h_m\} } \, 
    {\bf 1}_{\{ \overline{S}_i - S_i < a_i^{(S)}, \; \forall 0\le i\le n\} } \,
    {\bf 1}_{\{ \max_{1\le \ell < m} \Delta S_{\mathscrHS_{h_\ell}} \le r^\theta\} } 
    \Big] .
    \label{Zr:second_moment-3}
\end{eqnarray}

We proceed with (\ref{Zr:second_moment-2}). Recall from (\ref{aS}) that $a_i^{(S)} := \lambda_\ell$ if $\mathscrHS_{h_{\ell-1}} \le i< \mathscrHS_{h_\ell}$. In particular, $a_n^{(S)} = \lambda_m$ on the event $\{ h_{m-1} \le \overline{S}_n < h_m\}$, so $\ee^{\overline{S}_n -S_n} \le \ee^{\lambda_m}$ on $\{ h_{m-1} \le \overline{S}_n < h_m\} \cap \{\overline{S}_n - S_n < a_n^{(S)}\}$. Consequently,
\begin{eqnarray*}
 &&\E[E_\omega  (Z_r^2)]
    \le
    \ee^{-(1+o(1))(2r)^{1/2}}
    +
    2\ee^{2\varepsilon_1 \, r^{1/2}}\, 
    \sum_{n=0}^\infty \sum_{m=1}^k 
    \ee^{\lambda_m+2\varepsilon r^{1/2}- (2+o(1)) \sum_{\ell=m+2}^k \frac{r^\chi}{\lambda_\ell} } 
    \times
    \\
 && \times 
    \Q \Big(\{ h_{m-1} \le \overline{S}_n < h_m\}
    \cap
    \{ \overline{S}_i - S_i < a_i^{(S)}, \; \forall 0\le i\le n\}
    \cap
    \{ \max_{1\le \ell < m} \Delta S_{\mathscrHS_{h_\ell}} \le r^\theta\}
    \Big) .
\end{eqnarray*}

\noindent According to Claim \ref{claim:sum_j}, this yields that
\begin{eqnarray*}
    \E[E_\omega  (Z_r^2)]
 &\le& \ee^{-(1+o(1))(2r)^{1/2}} +
    \\
 && +
    2\ee^{2\varepsilon_1 \, r^{1/2}}\, 
    \sum_{m=1}^k 
    \ee^{\lambda_m+2\varepsilon r^{1/2}- (2+o(1)) \sum_{\ell=m+2}^k \frac{r^\chi}{\lambda_\ell} } 
    \times
    c_{14}r \, \ee^{- (1+o(1)) \sum_{\ell=1}^{m-1} \frac{r^\chi}{\lambda_\ell} } \, .
\end{eqnarray*}

\noindent By definition, $k := \lfloor r^{1-\chi}\rfloor$ and $\lambda_m:= (2r)^{1/2} \, (\frac{k-m+1}{k})^{1/2}$. Hence
$$
\lambda_m
- 
2\sum_{\ell=m+2}^k \frac{r^\chi}{\lambda_\ell}
-
\sum_{\ell=1}^{m-1} \frac{r^\chi}{\lambda_\ell} 
\sim
-(2r)^{1/2} \, .
$$

\noindent This completes the proof of the inequality (\ref{E(Z2)}) in Lemma \ref{l:Z}.\hfill$\Box$

\subsection{Proof of Lemma \ref{l:Z}: the inequality (\ref{(EZ)2})}
\label{subs:proof-second-lemma-(iii)}

We recall from (\ref{Zr:second_moment-3}) that
\begin{eqnarray*}
 &&\E[(E_\omega Z_r)^2]
    \le
    \ee^{-r}
    +
    \sum_{n=0}^\infty \sum_{m=1}^k 
    \ee^{2\varepsilon r^{1/2}- (2+o(1)) \sum_{\ell=m+2}^k \frac{r^\chi}{\lambda_\ell} } \,
    \E_{\Q} \Big[ 
    \ee^{-S_n} \,
    {\bf 1}_{\{ \min_{0\le i\le n} S_i \ge -\beta\} } \, \times
    \\
 &&\qquad\qquad \times 
    {\bf 1}_{\{ h_{m-1} \le \overline{S}_n < h_m\} } \, 
    {\bf 1}_{\{ \overline{S}_i - S_i < a_i^{(S)}, \; \forall 0\le i\le n\} } \,
    {\bf 1}_{\{ \max_{1\le \ell < m} \Delta S_{\mathscrHS_{h_\ell}} \le r^\theta\} } 
    \Big] .
\end{eqnarray*}

\noindent On the right-hand side, we throw away ${\bf 1}_{\{ \max_{1\le \ell < m} \Delta S_{\mathscrHS_{h_\ell}} \le r^\theta\} }$ by saying that it is bounded by 1. On the event $\{ h_{m-1} \le \overline{S}_n < h_m\}$, we have $a_n^{(S)} = \lambda_m$, so ${\bf 1}_{\{ \overline{S}_i - S_i < a_i^{(S)}, \; \forall 0\le i\le n\} } \le {\bf 1}_{\{ \overline{S}_n - S_n < \lambda_m\} }$. This leads to:
\begin{eqnarray}
    \E[(E_\omega Z_r)^2]
 &\le&\ee^{-r}
    +
    \sum_{m=1}^k 
    \ee^{2\varepsilon r^{1/2}- (2+o(1)) \sum_{\ell=m+2}^k \frac{r^\chi}{\lambda_\ell} } \,
    \E_{\Q} \Big[ \sum_{n=0}^\infty
    \ee^{-S_n} \,
    {\bf 1}_{\{ \min_{0\le i\le n} S_i \ge -\beta\} } \, \times
    \nonumber
    \\
 &&\qquad\qquad \times 
    {\bf 1}_{\{ h_{m-1} \le \overline{S}_n < h_m\} } \, 
    {\bf 1}_{\{ \overline{S}_n - S_n < \lambda_m\} } \,
    \Big] 
    \nonumber
    \\
 &=:& \ee^{-r} + \sum_{m=1}^k \Sigma_4^{(m)} \, ,
    \label{second-moment-5}
\end{eqnarray}

\noindent with obvious notation.

Fix $0<\varepsilon_5<1$. We use different estimates for $\Sigma_4^{(m)}$ on the right-hand side, depending on whether $m \le \lceil \varepsilon_5 k\rceil$ or not. 

{\bf First case: $1\le m\le \lceil \varepsilon_5\, k\rceil$.} In this case, we simply use ${\bf 1}_{\{ h_{m-1} \le \overline{S}_n < h_m\} } \le 1$ and ${\bf 1}_{\{ \overline{S}_n - S_n < \lambda_m\} } \le 1$, to see that for large $r$,
$$
\Sigma_4^{(m)}
\le
\ee^{2\varepsilon r^{1/2}- (2+o(1)) \sum_{\ell=m+2}^k \frac{r^\chi}{\lambda_\ell} } \, \E_{\Q} \Big[ \sum_{n=0}^\infty
\ee^{-S_n} \,
{\bf 1}_{\{ \min_{0\le i\le n} S_i \ge -\beta\} } \Big] .
$$

\noindent According to Lemma B.2 of A\"\i d\'ekon~\cite{elie-min}, for any $b>0$, there exists a constant $c_{20}(b)>0$, whose value depends also on $\beta$, such that
\begin{equation}
    \E_{\Q} \Big[ \sum_{j=1}^\infty \ee^{- b\, S_j}\, {\bf 1}_{\{ S_i\ge -\beta, \, \forall i \le j\} } \Big]
    \le 
    c_{20}(b) \, .
    \label{elie_appendice}
\end{equation}

\noindent Consequently, for all sufficiently large $r$,
$$
\Sigma_4^{(m)}
\le
c_{20}(1) \, \ee^{2\varepsilon r^{1/2}- (2+o(1)) \sum_{\ell=m+2}^k \frac{r^\chi}{\lambda_\ell} } \, .
$$

\noindent By (\ref{sum_over_lambda_j}) and (\ref{sum(hm-hm-1)/lambda_m=sqrt(2r)}), for $1\le m\le \lceil \varepsilon_5 k\rceil$, we have 
$$
\sum_{\ell=m+2}^k \frac{r^\chi}{\lambda_\ell}
=
\sum_{\ell=1}^k \frac{r^\chi}{\lambda_\ell}
-
\sum_{\ell=1}^{m+1} \frac{r^\chi}{\lambda_\ell}
=
(1+o(1)) (2r)^{1/2}
-
(2r^\chi)^{1/2}[k^{1/2} - (k-\lceil \varepsilon_5 k\rceil)^{1/2}] ,
$$

\noindent which is $(1+o(1))(1-\varepsilon_5)^{1/2}(2r)^{1/2}$, $r\to \infty$. Therefore,
\begin{equation}
    \sum_{m=1}^{\lceil \varepsilon_5\, k\rceil} \Sigma_4^{(m)} 
    \le
    c_{20}(1) \, \lceil \varepsilon_5\, k\rceil \, \ee^{2\varepsilon r^{1/2}- (2+o(1)) (1-\varepsilon_5)^{1/2}(2r)^{1/2} } \, .
    \label{second-moment-6}
\end{equation}

{\bf Second (and last) case: $\lceil \varepsilon_5\, k\rceil < m \le k$.} Since $m> \lceil \varepsilon_5\, k\rceil$, we have $h_{m-1} = (m-1)\, \frac{r}{k} \ge \varepsilon_5 r$. So on the event $\{ h_{m-1} \le \overline{S}_n < h_m\} \cap \{ \overline{S}_n - S_n < \lambda_m\}$, we have $S_n > \overline{S}_n - \lambda_m \ge h_{m-1} - \lambda_m \ge \varepsilon_5 r - \lambda_m$, which is greater than or equal to $\varepsilon_5 r - \lambda_1 = \varepsilon_5 r - (2r)^{1/2}$. Accordingly,
\begin{eqnarray*}
    \Sigma_4^{(m)}
 &\le&\ee^{2\varepsilon r^{1/2}- (2+o(1)) \sum_{\ell=m+2}^k \frac{r^\chi}{\lambda_\ell} } \, \E_{\Q} \Big[ \sum_{n=0}^\infty
    \ee^{-\frac12 S_n} \, \ee^{-\frac12 [\varepsilon_5 r - (2r)^{1/2}]} \,
    {\bf 1}_{\{ \min_{0\le i\le n} S_i \ge -\beta\} } \Big]
    \\
 &\le&\ee^{2\varepsilon r^{1/2}} \, \E_{\Q} \Big[ \sum_{n=0}^\infty
    \ee^{-\frac12 S_n} \, \ee^{-\frac12 [\varepsilon_5 r - (2r)^{1/2}]} \,
    {\bf 1}_{\{ \min_{0\le i\le n} S_i \ge -\beta\} } \Big] 
    \\
 &=& \ee^{2\varepsilon r^{1/2}-\frac12 [\varepsilon_5 r - (2r)^{1/2}]} \, \E_{\Q} \Big[ \sum_{n=0}^\infty
    \ee^{-\frac12 S_n} \, 
    {\bf 1}_{\{ \min_{0\le i\le n} S_i \ge -\beta\} } \Big] .   
\end{eqnarray*}

\noindent So by (\ref{elie_appendice}), we have $\Sigma_4^{(m)} \le c_{20}(\frac12) \, \ee^{2\varepsilon r^{1/2}-\frac12 [\varepsilon_5 r - (2r)^{1/2}]}$ for $\lceil \varepsilon_5\, k\rceil < m \le k$. As a consequence,
\begin{equation}
    \sum_{m= \lceil \varepsilon_5 k\rceil +1}^k \Sigma_4^{(m)}
    \le
    c_{20}(1/2)k \, \ee^{2\varepsilon r^{1/2}-\frac12 [\varepsilon_5 r - (2r)^{1/2}]} \, .
    \label{second-moment-7}
\end{equation}

Since $\E[(E_\omega Z_r)^2] \le \ee^{-r} + \sum_{m=1}^k \Sigma_4^{(m)}$ (see (\ref{second-moment-5})), it follows from (\ref{second-moment-6}) and (\ref{second-moment-7}) that 
\begin{eqnarray*}
    \E[(E_\omega Z_r)^2] 
 &\le& \ee^{-r}
    +
    c_{20}(1) \, \lceil \varepsilon_5\, k\rceil \, \ee^{2\varepsilon r^{1/2}- (2+o(1)) (1-\varepsilon_5)^{1/2}(2r)^{1/2} } +
    \\
 &&\qquad\qquad +
    c_{20}(1/2)k \, \ee^{2\varepsilon r^{1/2}-\frac12 [\varepsilon_5 r - (2r)^{1/2}]} \, .
\end{eqnarray*}

\noindent Recall that $k:= \lfloor r^{1-\chi}\rfloor$. Since $\varepsilon_5>0$ can be as close to $0$ as possible, this yields (\ref{(EZ)2}), and completes the proof of Lemma \ref{l:Z}.\hfill$\Box$

\appendix
\section{Probability estimates for one-di\-men\-sional random walks}
   \label{s:appendix}

Let $(\Omega, \, \mathscr{F}, \, \p)$ be a probability space. Let $S_0:=0$ and let $(S_i-S_{i-1}, \, i\ge 1)$ be a sequence of i.i.d.\ real-valued random variables defined on $(\Omega, \, \mathscr{F}, \, \p)$ with $\e(S_1)=0$ and $\sigma^2:= \e(S_1^2) \in (0, \, \infty)$. We write
$$
\overline{S}_j
:=
\max_{0\le i\le j} S_i,
\qquad j\ge 0.
$$

\noindent For any $b\in \r$, let\footnote{For $b>0$, $\HH_b$ is nothing else but $\mathscrHS_b$ defined in (\ref{mathscrHS_r}).}
$$
\HH_b
:=
\inf\{ i\ge 1: \, S_i \ge b\},
\qquad
\HH^-_{b}
:=
\inf\{ i\ge 1: \, S_i \le b\} \, .
$$

\noindent Applying (2.6) of Borovkov and Foss~\cite{borovkov-foss} to the ladder heights, we immediately see that the assumption $\e(S_1^2)<\infty$ ensures that $\e(S_{\HH_b})<\infty$ for all $b\ge 0$, and that there exists a constant $c_{21}>0$ satisfying $\e(S_{\HH_b}-b) \le c_{21} (b+1)$ for all $b\ge 0$.

\begin{lemma}
\label{l:RW-exit-place}

 {\rm (i)} Assume $\e(|S_1|^3)<\infty$. 
 There exists a constant $c_{22}>0$ such that for any $a\ge 0$ and $b\ge 0$ with $a+b>0$,
 \begin{equation}
     \frac{b-c_{22}}{a+b}
     \le
     \p\{\HH_a < \HH^-_{-b}\}
     \le
     \frac{b+c_{22}}{a+b} \, .
     \label{RW-exit-place}
 \end{equation}

 {\rm (ii)} Assume $\e(|S_1|^{3+\delta})<\infty$ for some $\delta>0$. Then for any $a\ge 0$,
 \begin{equation}
     \p\{\HH^-_{-b} < \HH_a\}
     \; \sim \;
     \frac{\e(S_{\HH_a})}{b} \, ,
     \qquad b \to \infty.
     \label{RW-exit-tail}
 \end{equation}

\end{lemma}

\noindent {\it Proof.} We follow the same argument as in \cite{AHZ}.

(i) Since $\e(|S_1|^3)<\infty$, it is known (Mogulskii~\cite{mogulskii73}) that $\sup_{b>0} \e(-S_{\HH^-_{-b}}-b) <\infty$. 

By the optional stopping theorem, $0= \e(S_{\HH_a \wedge \HH^-_{-b}}) = \e[(S_{\HH_a}-S_{\HH^-_{-b}}) \, {\bf 1}_{\{ \HH_a < \HH^-_{-b} \} }] + \e(S_{\HH^-_{-b}}) \ge (a+b)\, \p\{ \HH_a < \HH^-_{-b} \} -b - \e(- S_{\HH^-_{-b}}-b) \ge (a+b)\, \p\{ \HH_a < \HH^-_{-b} \} -b - c_{23}$ where $c_{23} := \sup_{b>0} \e(-S_{\HH^-_{-b}}-b) <\infty$. This yields the second inequality in (\ref{RW-exit-place}). Considering $(-S_n)$ in place of $(S_n)$ (and exchanging the roles of $a$ and $b$) yields the first inequality.

(ii) Again, by the optional stopping theorem, $0= \e(S_{\HH_a \wedge \HH^-_{-b}}) = -b \, \p \{\HH^-_{-b} < \HH_a\} + \e(S_{\HH_a}) + \e\{ [(S_{\HH^-_{-b}}+b)-S_{\HH_a}] \, {\bf 1}_{\{ \HH^-_{-b} < \HH_a \} } \}$, which leads to
\begin{equation}
    b \, \p \{\HH^-_{-b} < \HH_a\}
    =
    \e(S_{\HH_a})
    +
    \e\{ [\, |S_{\HH^-_{-b}}+b| + S_{\HH_a}] \, {\bf 1}_{\{\HH^-_{-b} < \HH_a\} } \} \, .
    \label{RW-exit-tail-prepa}
\end{equation}

\noindent We let $b\to \infty$. We have $\p \{\HH^-_{-b} < \HH_a\} \to 0$ (by (\ref{RW-exit-place})), whereas $\sup_{b>0} \e (\, | S_{\HH^-_{-b}} + b|^{1+\delta})<\infty$ and $\e[(S_{\HH_a})^{1+\delta})<\infty$ (which is a consequence of the assumption $\e(|S_1|^{3+\delta})<\infty$; see Mogulskii~\cite{mogulskii73}). By H\"older's inequality, $\e\{ [\, |S_{\HH^-_{-b}}+b| + S_{\HH_a}] \, {\bf 1}_{\{\HH^-_{-b} < \HH_a\} } \} \to 0$. So (\ref{RW-exit-tail-prepa}) implies (\ref{RW-exit-tail}).\hfill$\Box$

\begin{lemma}
\label{l:RW-chute-lb}
 
 Assume $\e(|S_1|^3)<\infty$.
 There exist constants $c_{24}>0$, $c_{25}>0$ and $c_{26}>0$ such that for all $r\ge 1$ and $\lambda\ge c_{24}$, we have
 \begin{equation}
     \p\Big\{ \overline{S}_j - S_j < \lambda, \, S_j \ge 0, \, \forall 0\le j\le \HH_r\Big\}
     \ge
     c_{25} \, \exp\Big( - \frac{r}{\lambda} - \frac{c_{26}\, r}{\lambda^{3/2}} \Big).
     \label{RW-chute-lb}
 \end{equation}

\end{lemma}

\noindent {\it Proof.} Let $c_{22}>0$ be the constant in Lemma \ref{l:RW-exit-place}. Since $\e(S_1)=0$ and $\e(S_1^2)>0$, there exist $c_{27}>0$ and $c_{28}\in (0, \, 1)$ such that $\p \{ S_1 \ge c_{27}\} \ge c_{28}$, so that
$$
\p\{ \HH_{c_{22}+1} < \HH^-_0\}
\ge
\p \Big\{ S_i - S_{i-1} \ge c_{27}, \; \forall 1\le i\le \lceil \frac{c_{22}+1}{c_{27}} \rceil \Big\}
\ge 
c_{28}^{\lceil \frac{c_{22}+1}{c_{27}} \rceil} 
=:
c_{29} > 0.
$$

Let $y>0$ and let $r_k := (c_{22}+1)+ yk$, for $0\le k\le N := \lceil \frac{r}{y} \rceil$.

Let $E_{(\ref{RW-chute-lb})}:= \{ \overline{S}_j - S_j < \lambda, \, S_j \ge 0, \, \forall 0\le j\le \HH_r \}$. Since $r_N \ge r$, $E_{(\ref{RW-chute-lb})}$ will be realized if $\HH_{r_0} < \HH^-_0$ and if for all $0\le k\le N-1$, the following is true: after hitting $[r_k, \, \infty)$ for the first time, the walk $(S_n)$ hits $[r_{k+1}, \, \infty)$ before hitting $(-\infty, \, r_k-\lambda]$. Applying the strong Markov property gives that ($\p_x$ being the probability under which the random walk starts at $x$; so $\p_0=\p$)
$$
\p(E_{(\ref{RW-chute-lb})})
\ge
\p\{ \HH_{r_0} < \HH^-_0 \} \times \prod_{k=0}^{N-1} \p_{r_k} \{ \HH_{r_{k+1}} < \HH^-_{r_k - \lambda} \}
\ge
c_{29} \prod_{k=0}^{N-1} \p_{r_k} \{ \HH_{r_{k+1}} < \HH^-_{r_k - \lambda} \} \, .
$$

\noindent [We do not need to worry about overshoots, because $x\mapsto \p_x \{ \HH_{r_{k+1}} < \HH^-_{r_k - \lambda} \}$ is non-decreasing for $x\in [r_k, \, \infty)$.] 

Since $\p_{r_k} \{ \HH_{r_{k+1}} < \HH^-_{r_k - \lambda} \} = \p \{ \HH_{r_{k+1}-r_k} < \HH^-_{- \lambda} \} = \p \{ \HH_y < \HH^-_{- \lambda} \}$, it follows from Lemma \ref{l:RW-exit-place} that (with $\lambda$ sufficiently large such that $\lambda > y+c_{22}$)
$$
\p_{r_k} \{ \HH_{r_{k+1}} < \HH^-_{r_k - \lambda} \}
\ge
\frac{\lambda - c_{22}}{y+\lambda}
=
1- \frac{y+c_{22}}{y+\lambda}
\ge
1- \frac{y+c_{22}}{\lambda} ,
$$

\noindent which is greater than or equal to $\exp[ -\frac{y+c_{22}}{\lambda} - (\frac{y+c_{22}}{\lambda})^2 ]$ if $\frac{y+c_{22}}{\lambda} \le \frac12$ (by the elementary inequality that $1-x \ge \ee^{-x-x^2}$ for $0\le x\le \frac12$). Since $N \le \frac{r}{y} +1 = \frac{r+y}{y}$, we obtain:
$$
\p(E_{(\ref{RW-chute-lb})})
\ge
c_{29} \exp \Big[ - \frac{y+c_{22}}{\lambda}\frac{r+1}{y} - \frac{(y+c_{22})^2}{\lambda^2}\frac{r+1}{y} \Big] .
$$

\noindent We choose $\lambda\ge 1$ and $r\ge 1$. We note that $\frac{y+c_{22}}{\lambda}\frac{r+1}{y} = \frac{r}{\lambda} + \frac{1}{\lambda} + \frac{c_{22}}{\lambda}\frac{r+1}{y} \le \frac{r}{\lambda} + 1 + \frac{2c_{22} r}{\lambda y}$, and that if $y\ge c_{22}$, $\frac{(y+c_{22})^2}{\lambda^2}\frac{r+1}{y} \le \frac{4y^2}{\lambda^2}\frac{2r}{y} = \frac{8ry}{\lambda^2}$. So, taking $y:= \lambda^{1/2}$ yields that
$$
\p(E_{(\ref{RW-chute-lb})})
\ge
c_{29} \exp \Big[ - \frac{r}{\lambda} - 1  - \frac{2c_{22} r}{\lambda^{3/2}} - \frac{8r}{\lambda^{3/2}} \Big] ,
$$

\noindent proving the lemma.\hfill$\Box$

\medskip

The next lemma says that, under sufficient integrability conditions, the main term $\frac{r}{\lambda}$ within the exponential function in Lemma \ref{l:RW-chute-lb} is, in some sense, optimal:

\begin{lemma}
\label{l:RW-chute-ub}
 
 Assume that $\e(\ee^{\delta S_1})<\infty$ for some $\delta>0$.
 For any $\varepsilon>0$, there exist constants $c_{30}>0$ and $c_{31}>0$ such that for all $r\ge 1$ and $\lambda\ge c_{30}$, we have
 \begin{equation}
     \p\Big\{ \overline{S}_j - S_j < \lambda, \, \forall 0\le j\le \HH_r\Big\}
     \le
     c_{31} \, \exp\Big( - (1-\varepsilon) \frac{r}{\lambda} \Big).
     \label{RW-chute-ub}
 \end{equation}

\end{lemma}

\noindent {\it Proof.} Let $\tau_0 := 0$ and for any $k\ge 1$, let $\tau_k := \inf\{ i> \tau_{k-1}: \, S_i \ge S_{\tau_{k-1}}\}$ be the $k$-th ascending ladder epoch. Let $\p_{(\ref{RW-chute-ub})}$ denote the probability expression on the left-hand side of (\ref{RW-chute-ub}). For any $k\ge 1$, we have
$$
\p_{(\ref{RW-chute-ub})}
\le
\p\{ S_{\tau_k} \ge r\}
+
\p\Big\{ S_{\tau_{i-1}} - \min_{\tau_{i-1} \le j\le \tau_i} S_j < \lambda, \; \forall 1\le i\le k\Big\} .
$$

\noindent We now estimate the two probability expressions on the right-hand side.

For the first probability expression, we write $S_{\tau_k} = \sum_{i=1}^k (S_{\tau_i} - S_{\tau_{i-1}})$, and observe that $(S_{\tau_i} - S_{\tau_{i-1}}, \, i\ge 1)$ is a sequence of i.i.d.\ random variables, with $\e(\ee^{a S_{\tau_1}})<\infty$ for all $a< \delta$. So we take 
$$
k
= 
k(r, \, \varepsilon) 
:= 
\Big\lceil \, \frac{1-\varepsilon}{\e(S_{\tau_1})} \, r \, \Big\rceil \, ;
$$

\noindent there exist constants $c_{32}>0$ and $c_{33}>0$, depending on $\varepsilon$, such that $\p\{ S_{\tau_{k(r, \, \varepsilon)}} \ge r\} \le c_{32} \, \ee^{-c_{33} \, r}$ for all $r\ge 1$. 

For the second probability expression (now with $k:= k(r, \, \varepsilon)$), we use the fact that $(S_{\tau_{i-1}} - \min_{\tau_{i-1} \le j\le \tau_i} S_j, \, i\ge 1)$ is also a sequence of i.i.d.\ random variables, having the same distribution as $-\min_{0\le j\le \tau_1}S_j$; accordingly,
$$
\p\Big\{ S_{\tau_{i-1}} - \min_{\tau_{i-1} \le j\le \tau_i} S_j < \lambda, \; \forall 1\le i\le k(r, \, \varepsilon)\Big\}
=
\Big[ \, \p\Big\{ -\min_{0\le j\le \tau_1}S_j < \lambda \Big\} \Big]^{k(r, \, \varepsilon)} \, .
$$

\noindent Since $\tau_1 = \HH_0$ and $\{ -\min_{0\le j\le \tau_1}S_j < \lambda \} = \{ \HH_0 < \HH_{-\lambda}^-\}$, we are entitled to apply (\ref{RW-exit-tail}) to see that for all sufficiently large $\lambda$ (say $\lambda \ge \lambda_0$), $\p \{ -\min_{0\le j\le \tau_1}S_j < \lambda \} \le 1- (1-\varepsilon) \frac{\e(S_{\tau_1})}{\lambda}$. Hence for $\lambda\ge \lambda_0$,
$$
\p\Big\{ S_{\tau_{i-1}} - \min_{\tau_{i-1} \le j\le \tau_i} S_j < \lambda, \; \forall 1\le i\le k(r, \, \varepsilon)\Big\}
\le
\Big( 1- (1-\varepsilon) \frac{\e(S_{\tau_1})}{\lambda} \Big)^{k(r, \, \varepsilon)} \, ,
$$

\noindent which is bounded by $\exp [ - (1-\varepsilon) \frac{\e(S_{\tau_1})}{\lambda}\, k(r, \, \varepsilon)]$. Assembling these pieces yields that for $r\ge 1$ and $\lambda\ge \lambda_0$,
$$
\p_{(\ref{RW-chute-ub})}
\le
c_{32} \, \ee^{-c_{33} \, r}
+
\exp \Big[ - (1-\varepsilon) \frac{\e(S_{\tau_1})}{\lambda}\, k(r, \, \varepsilon) \Big] \, ,
$$

\noindent which yields (\ref{RW-chute-ub}) as $\varepsilon>0$ is arbitrary.\hfill$\Box$

\section{Proof of \eqref{cvVXn}}
\label{s:cvVXn}

The proof of \eqref{cvVXn} relies on several results from \cite{yzlocaltree}, and we use the notation therein. Let $t>0$. By \cite[Theorems 2.1 and 2.7]{yzlocaltree}, $$P_\omega\Big(\frac{ V(X_n)}{\log n} \le t \Big)= \frac{\sigma^2}{2 D_\infty \log n} \sum_{ x \in \T} \ee^{-V(x)} {\bf 1}_{\{   V(x) \le t \log n, \, x < {\mathscr L}_n \}}+ o_{\P^*}(1),$$


\noindent where $x < {\mathscr L}_n$ means that for all $y\in \, ]\!] \varnothing,\, x]\!]$, $\sum_{z\in \, ]\!] \varnothing,\, y]\!]} \ee^{V(z)- V(y)} \le   n$,  and  $o_{\P^*}(1)$ denotes a quantity which converges to $0$ in $\P^*$-probability as $n\to\infty$. Let for any $s, \lambda  >0$, $$W_n^{(s, \lambda)}:= \sum_{|x|=n} \ee^{-V(x)} {\bf 1}_{\{V(x) \le s, \,    \max_{y\in \, ]\!] \varnothing,\, x]\!]} (\overline V(y)- V(y)) \le \lambda\}}.$$

\noindent A line-by-line analogue of the proof of \cite[Lemma 4.1]{yzlocaltree} yields that  $$ \lim_{\lambda\to\infty} \frac1\lambda \sum_{k=1}^\infty W_k^{(t \lambda, \lambda)} =\left(\frac{8}{\pi}\right)^{1/2} \, \frac{D_\infty}{\sigma^2}\,  \E\left( \min\left( \frac1{\tt m_1^{\#}}, \frac{t}{\tt m_1}\right)\right), \qquad \mbox{in $\P^*$-probability}.$$

From this point, we can closely follow the step-by-step arguments in the proof of Corollary 2.3 of \cite{yzlocaltree}. Let us give an outline.  First, by using \cite[Equations (4.6) and (4.7)]{yzlocaltree}, we have that for $B>b>0$, \begin{align*}   & \frac{\sigma^2}{2 D_\infty \log n} \sum_{ x \in \T} \ee^{-V(x)} {\bf 1}_{\{   V(x) \le t \log n, \, x < {\mathscr L}_n \}}
\\
=& \frac{\sigma^2}{2 D_\infty \log n} \sum_{b (\log n)^2 \le |x| \le B (\log n)^2} \ee^{-V(x)} {\bf 1}_{\{   V(x) \le t \log n, \, x< {\mathscr L}_n\}}+ o_{\P^*, b, B}(1),
\end{align*}

\noindent where the term $o_{\P^*, b, B}(1)$ denotes a quantity which converges to $0$ in $\P^*$-probability when first $n\to\infty$, then $b \to 0$ and $B\to \infty$. Exactly as in the proof of Corollary  2.3 of  \cite{yzlocaltree},  we can use the following inequalities 
\begin{align*}   
\sum_{b (\log n)^2 \le |x| \le B (\log n)^2} \ee^{-V(x)} {\bf 1}_{\{   V(x) \le t \log n, \, x< {\mathscr L}_n\}}
&\le \sum_{k=b (\log n)^2 }^{ B (\log n)^2} W_k^{(t \log n,  \log n)}, 
\\
\sum_{b (\log n)^2 \le |x| \le B (\log n)^2} \ee^{-V(x)} {\bf 1}_{\{   V(x) \le t \log n, \, x< {\mathscr L}_n\}}
&\ge \sum_{k=b (\log n)^2 }^{ B (\log n)^2} W_k^{( t \log n, \log (n/B(\log n)^2))}, \end{align*}

\noindent to deduce that $$ \frac{\sigma^2}{2 D_\infty \log n} \sum_{ x \in \T} \ee^{-V(x)} {\bf 1}_{\{   V(x) \le t \log n, \, x < {\mathscr L}_n \}}
= \left(\frac{2}{\pi}\right)^{1/2} \, \E\left( \min\left( \frac1{\tt m_1^{\#}}, \frac{t}{\tt m_1}\right)\right)+o_{\P^*}(1). $$
 
 \noindent Then \eqref{cvVXn} follows. \hfill$\Box$

\end{document}